%September 19, 2014
\input amstex
\input psfig.sty
\magnification 1200
\TagsOnRight
\def\qed{\ifhmode\unskip\nobreak\fi\ifmmode\ifinner\else
 \hskip5pt\fi\fi\hbox{\hskip5pt\vrule width4pt
 height6pt depth1.5pt\hskip1pt}}
\NoBlackBoxes
\baselineskip 20 pt
\parskip 5 pt
\def\stretch {\noalign{\medskip}}
\define \bC {\bold C}
\define \bCp {\bold C^+}
\define \bCm {\bold C^-}
\define \bCpb {\overline{\bold C^+}}
\define \bCmb {\overline{\bold C^-}}
\define \ds {\displaystyle}
\define \bR {\bold R}
\define \bm {\bmatrix}
\define \endbm {\endbmatrix}

{
\parskip 5 pt

\centerline {\bf INVERSE PROBLEMS FOR SELFADJOINT SCHR\"ODINGER OPERATORS}
\centerline {\bf ON THE HALF LINE WITH COMPACTLY-SUPPORTED POTENTIALS}

\baselineskip 12 pt
\vskip 5 pt
\centerline {Tuncay Aktosun}
\centerline {Department of Mathematics}
\centerline {University of Texas at Arlington}
\centerline {Arlington, TX 76019-0408, USA}
\centerline {aktosun\@uta.edu}

\vskip 2 pt

\centerline {Paul Sacks}
\centerline {Department of Mathematics}
\centerline {Iowa State University}
\centerline {Ames, IA 24061, USA}
\centerline {psacks\@iastate.edu}

\vskip 2 pt

\centerline {Mehmet Unlu}
\centerline {Department of Mathematics}
\centerline {University of Texas at Arlington}
\centerline {Arlington, TX 76019-0408, USA}
\centerline {mehmet.unlu\@mavs.uta.edu}

%\vskip 2 pt

\noindent {\bf Abstract}:
For a selfadjoint Schr\"odinger operator
on the half line with a real-valued,
integrable, and compactly-supported potential, it is investigated whether
the boundary parameter at the origin and the potential can
uniquely be determined by the scattering matrix or
by the absolute
value of the Jost function known at positive energies,
without having the bound-state
information.
It is proved that, except in one special case where the scattering
 matrix has no bound states and its value is $+1$ at zero energy,
 the determination by the scattering matrix is unique. In the special
 case, it is shown that there are exactly two distinct sets consisting
 of a potential and a boundary
 parameter yielding the same scattering matrix, and a characterization of
 the nonuniqueness is provided.
 A reconstruction from the scattering matrix is outlined
 yielding all the corresponding potentials and boundary parameters.
The concept of ``eligible resonances" is introduced, and such resonances
correspond to real-energy resonances that can be converted into
bound states via a Darboux
transformation without changing the compact support of the potential.
It is proved that the determination of the boundary parameter and the potential by
the absolute
value of the Jost function is unique up to the inclusion of eligible
resonances. Several equivalent characterizations are
provided to determine whether a resonance is eligible or ineligible.
A reconstruction from the absolute value of the Jost function
is given, yielding all the corresponding potentials and boundary
parameters.
The results obtained are illustrated with various explicit examples.


\par\noindent {\bf Short title:} Inverse problem with compactly-supported potentials
%for the Schr\"odinger equation

}
\newpage

\noindent {\bf 1. INTRODUCTION}
\vskip 3 pt

In this paper we consider the half-line Schr\"odinger operator
with the general selfadjoint boundary condition at the origin
when the potential is real valued, integrable, and compactly supported. We
examine the inverse problem of recovery of the potential and
boundary condition from two distinct types of input data,
investigate whether the determination from each input data set is unique,
present the characterization of the nonuniqueness
if the unique determination is not possible, and provide a procedure
to reconstruct all the corresponding potentials and boundary conditions from each input data set.

The first set of input data we use is the scattering matrix known
at all positive energies, but without any
explicit information on the bound states. The second input data set we use
is the absolute value of the so-called Jost function
given at all positive energies, but again without
any explicit information on the bound states. Assuming that the existence problem
is solved, i.e. by assuming that there exists at least one set consisting of a potential and
a boundary condition corresponding to our input data, we investigate
whether we have two or more distinct sets containing a potential and a boundary
condition corresponding to our input data and provide a reconstruction of
all such sets.

Our inverse
scattering problem can be paraphrased as follows:
To what extent, can the lack of bound-state information in our
input data set be compensated by the knowledge that the potential
is compactly supported? We certainly need to restrict our study to
a specific class of potentials so that the problem under study is
mathematically well stated. Real-valued, integrable potentials
naturally arise [1,7,8,14-16] in the theory of inverse
problems for Schr\"odinger operators
on the half line. The potentials of compact support
appear in our analysis because for such potentials
the corresponding
Jost function has an analytic extension from the real axis
to the entire complex plane. Such an analytic extension is crucial in
our analysis in order to compensate for
the lack of bound-state information in our data.

A motivation to study our inverse problems comes from
the inverse problem of determining the radius of the
human vocal tract from sound-pressure measurements at the lips [4].
The vocal tract radius as a function of the distance
 from the glottis is related to the potential of the Schr\"odinger
 equation, the length of the vocal tract corresponds to the length
 of the support interval of the potential, the behavior of the vocal
 tract at the glottis is accounted for by the selfadjoint boundary condition
 for the Schr\"odinger operator,
and the sound pressure at the lips as a
 function of the sound frequency  is related to
 the absolute value of
 the Jost function. The human speech consists of
 phonemes, and during the utterance of a phoneme if the upper lip opens
 downward (i.e. when the slope of the radius of the vocal tract at
 the upper lip is negative) as in the utterance of
 the vowel /o/, then the corresponding Schr\"odinger operator
 has one bound state, and the Schr\"odinger operator has no bound states
if the slope of the radius function at the upper lip is positive or zero
as in the utterance of /a/ or /u/, respectively.

There are two main methods to solve the inverse problem for a
selfadjoint Schr\"odinger operator
on the half line. The first is the Marchenko method [1,7-10,14,15], and
it uses the input data set consisting of the scattering matrix and
the bound-state information. In the Marchenko method the bound-state
information consists of the bound-state energies and the so-called
bound-state norming constants. The second method is the Gel'fand-Levitan method
[7,8,11,14,15], and that method uses the input data set consisting of
the absolute value of the Jost function and the bound-state
information. In the Gel'fand-Levitan method, the bound-state
information consists of the bound-state energies (such energies
are the same as the bound-state energies used in the Marchenko method)
and the bound-state norming constants (the Marchenko norming constants
and the Gel'fand-Levitan norming constants differ
 from each other even though they are related to each other).
In this paper, we consider the Marchenko recovery method when the
bound-state information is absent from the standard Marchenko input data but instead we know that
the corresponding potential is compactly supported. Similarly,
we consider the Gel'fand-Levitan method when the Gel'fand-Levitan input data
set does not contain the bound-state information but instead we know that
the corresponding potential is compactly supported.

The results proved in our paper are analogous to some results related the
full-line Schr\"odinger equation where the bound-state information is
missing from the input data. For example, a real-valued, integrable
potential with a finite first moment is uniquely determined [2,17] from
the corresponding left (right) reflection
coefficient alone if the support of the potential
is confined to the right (left) half line , or such a potential
is uniquely determined [3,13] by the data
consisting of the left (right) reflection
coefficient and knowledge of the potential on the left (right)
half line.

The analysis of the two inverse problems under study in our paper
turns out to have impact on other related
problems. One contribution of our study is in the area of
resonances for selfadjoint Schr\"odinger operators
on the half line. The nonzero zeros of the analytic extension of the
Jost function to the complex plane correspond to either bound states or
resonances. If such zeros are located
in the open upper-half complex plane, they correspond to bound states. It is known [1,7,8,14,15]
that each such bound-state zero is simple and that the number of such zeros is either zero
or a positive integer.
If the zeros of the Jost function
are located in the open lower-half complex plane, then those zeros
correspond to resonances. Equivalently, the poles of the
meromorphic extension of the scattering matrix correspond to bound states
if such poles occur in the open upper-half complex plane, and those poles of
the scattering matrix occurring
in the open lower-half complex plane correspond to resonances. The number of resonances
can be zero, one, or countably infinite.
%, which
%is a consequence of the fact that the Jost function is analytic
%in the entire complex plane.
A zero of the Jost function
corresponding to a resonance may or may not be simple.
The only real zero of the Jost function can occur at zero, and such a zero
is simple.

In our paper, we specifically deal with
resonances corresponding to the zeros of the Jost function
on the negative imaginary axis in the complex plane,
i.e. with real-energy resonances. In our analysis,
in a natural way, we are prompted to classify such resonances into two
mutually exclusive groups. The first group consists of
``eligible" resonances because such resonances can be converted into
bound states through a Darboux transformation [8,9,15] without changing the
compact support of the potential. The remaining resonances occurring on
the negative imaginary axis consist of ``ineligible" resonances because such
resonances cannot be converted into bound states under a Darboux transformation without changing the
compact support of the potential. It is remarkable that
ineligible resonances still remain ineligible if we add or remove
any number of bound states via a Darboux transformation without changing the
compact support of the potential. On the other hand, an eligible resonance
either remains eligible or is converted into a bound state if we add
any number of bound states via a Darboux transformation without changing the
compact support of the potential. Similarly, a bound state removed via a Darboux
transformation is converted into an eligible resonance.

Consider the sequence where each element in the sequence
consists of a potential and a boundary parameter in such a way that one element in the sequence is connected to another
 element through a number of Darboux transformations related to removing or
adding bound states without changing the compact support of the
 potentials.
For such a sequence, we define the ``maximal number of eligible resonances"
as the number of eligible resonances corresponding
to a pair with no bound states. Without causing any ambiguity, for
any term in the sequence
we can define the maximal number of eligible resonances
as the maximal number of
eligible resonances associated with the sequence itself.
Hence, for any term in the sequence
the sum of the number of eligible resonances and the number of bound states
must be equal to the maximal number of eligible resonances.
It turns out that each eligible resonance is simple
in the sense that the corresponding zero of the related Jost function
is a simple zero.
Hence, we do not need to be concerned about the multiplicity
of an eligible resonance. On the other hand, an ineligible
resonance does not need be simple, i.e. the corresponding zero
of the related Jost function
may not necessarily be a simple zero.

It is remarkable that
the identification of each
resonance on the negative imaginary axis either as eligible or ineligible arises
in a natural way
and is motivated by physics, and the identification
can be unambiguously given mathematically. One could certainly insist on converting an ineligible resonance
into a bound state, but in that case the resulting potential would no longer be
in the original class; either the compact support property would be lost or the
resulting potential would no longer be integrable. We illustrate the concepts of
eligible and ineligible resonances with some explicit examples in Section~6.

In the recovery of the potential and the selfadjoint boundary
condition from the scattering matrix $S_\theta(k),$ we summarize our main findings as follows. We have the unique
recovery, except in one special case. That special case occurs when
there are precisely two simultaneous constraints on $S_\theta(k),$
namely $S_\theta(0)=+1$ and at the same time there are
no bound-state poles associated with $S_\theta(k).$ The latter restriction is equivalent to the statement
that $S_\theta(k)$ has no poles
on the positive imaginary axis in the complex plane. In the special
case, it turns out that the scattering matrix
corresponds to exactly two distinct sets, each consisting of a potential and
a selfadjoint boundary condition.
Interestingly, when such a nonuniqueness occurs, the boundary condition in one set
must be the Dirichlet boundary condition and the boundary condition
in the other set must be a Neumann boundary condition. In Section~4
we further explore the nonuniqueness in the special case and provide
an interpretation of the nonuniqueness by viewing the compactly-supported
potential in the context of the corresponding full-line Schr\"odinger operator.
We then find that one of the nonunique potentials corresponds to the
reflection coefficient $R(k)$ and the other corresponds to
$-R(k),$ and this occurs when the corresponding full-line
Schr\"odinger operator has no bound states and is exceptional, i.e.
$R(0)\ne -1.$
In Section~6 we illustrate the nonuniqueness with an explicit example.

Concerning the recovery of the potential and the selfadjoint boundary
condition from the absolute value of the Jost function, we have the unique
recovery up to the inclusion of eligible resonances. From
our input data set
we are able to uniquely determine all eligible resonances. Let us use $M$ to denote
the maximal number of eligible resonances corresponding
to our input data set. We find that there are precisely $2^M$ distinct sets, each consisting of a potential and
a selfadjoint boundary condition, corresponding to the same input data.
We note [19] that $M$ can be infinite for our selfadjoint Schr\"odinger operator
on the half line when the potential is real valued, integrable, and compactly supported. A further minimal assumption [19] on the potential guarantees
that $M$ is finite. In Section~5 we present the details of the recovery from
the absolute value of the Jost function and elaborate on the $2^M$-fold nonuniqueness.

Our paper is organized as follows. In Section~2 we provide the preliminary
mathematical tools needed to analyze the two inverse problems under study. This is done by
introducing the half-line Schr\"odinger operator,
the selfadjoint boundary condition at the origin, the
Jost solution and the regular solution to the
half-line Schr\"odinger equation, the associated Jost function,
the scattering matrix, the bound
states, the norming constants, the resonances, and the relevant properties
of all such quantities.
In Section~3 we introduce the Darboux transformations to add or remove bound states,
obtain a few results related to the Darboux transformations
for potentials of compact support,
and provide several equivalent characterizations of eligible resonances.
In Section~4 we analyze the recovery
of the potential and the boundary condition from the scattering matrix alone.
We show that the recovery of the corresponding potential and
the boundary parameter is unique except in one special case, and we
characterize the double nonuniqueness in that special case.
In Section~5 we study the recovery
problem from the absolute value of the Jost function. We show that
the recovery is unique up to the inclusion of eligible resonances, which
is equivalent to having a $2^M$-fold nonuniqueness, with $M$
denoting the maximal number of eligible resonances. Finally, in Section~6 we
provide some explicit examples to illustrate the
theoretical results presented in Sections~3-5.

\vskip 10 pt
\noindent {\bf 2. PRELIMINARIES}
\vskip 3 pt

In this section we present the preliminaries
needed to prove the main results given
in Sections~3-5. We use $\bR$ to denote the real axis,
let $\bR^+:=(0,+\infty),$ use $\bC$ for the complex plane,
$\bCp$ for the open upper-half complex plane, $\bCm$ for the
open lower-half complex plane, $\bCpb:=\bCp\cup\bR,$ and $\bCmb:=\bCm\cup\bR.$

Consider the half-line Schr\"odinger equation
$$-\psi''+V(x)\,\psi=k^2\psi,\qquad x\in\bR^+,\tag 2.1$$
where the prime denotes the $x$-derivative and the potential $V$ is assumed
to belong to class $\Cal A$ defined as
$$\Cal A:=\left\{  V: \ V(x)\in\bR,\ V(x)\equiv 0 \text{ for }
x>b,\ \int_0^b dx\,|V(x)|<+\infty\right\} ,\tag 2.2$$
i.e. $V$ is real valued and integrable and it vanishes
when $x>b$ for some nonnegative $b.$
We obtain a selfadjoint Schr\"odinger operator on the half line
by supplementing (2.1) and (2.2) with the
general selfadjoint boundary condition at $x=0$ given by [7,11,14,15]
$$(\sin \theta)\,\psi'(0)+(\cos  \theta)\,\psi(0)=0,\tag 2.3$$
where the boundary parameter $\theta$ is a fixed real constant
 in the interval $(0,\pi].$ The
case $\theta=\pi$ in (2.3) corresponds to the Dirichlet boundary
condition $\psi(0)=0$ and
a case with $\theta\in(0,\pi)$
corresponds to a non-Dirichlet boundary condition
$\psi(0)\ne 0.$ The non-Dirichlet case
with $\theta=\pi/2$ in (2.3), i.e. $\psi'(0)=0,$ is known as the Neumann boundary condition.
The Dirichlet case arises especially when (2.1) is related to the
three-dimensional Schr\"odinger equation with a spherically symmetric potential.
On the other hand, there are various vibration problems [12] where a
non-Dirichlet boundary condition is more appropriate to use.
The non-Dirichlet case also arises in the inverse problem of
determining the shape of a human vocal tract from sound-pressure measurements
at the lips [4].

The so-called Jost solution associated with (2.1) and (2.2) is usually denoted by
$f(k,x),$ and it satisfies
$$f(k,x)=e^{ikx},\quad
f'(k,x)=ike^{ikx},\qquad x\ge b.\tag 2.4$$
For each fixed $x\in\bR^+\cup \{ 0\},$ the quantities
$f(k,x)$ and $f'(k,x)$ have analytic extensions [7-9,14,15] from $k\in\bR$ to
$k\in\bC$ as a consequence of $V$ belonging to class $\Cal A.$ Thus, for each fixed $x,$
the Jost function $f(k,x)$
has a Taylor series expansion around any $k$-value in $\bC.$

The so-called regular solution associated with (2.1)-(2.3), denoted
by $\varphi_\theta(k,x),$ satisfies the
initial conditions
$$\cases\varphi_\theta(k,0)=
1,\quad \varphi'_\theta(k,0)=-\cot\theta,&\qquad \theta\in (0,\pi),\\
\stretch
\varphi_\theta(k,0)=
0,\quad \varphi'_\theta(k,0)=1,&\qquad \theta=\pi.
\endcases\tag 2.5$$
The subscript $\theta$ in $\varphi_\theta(k,x)$ indicates the dependence on
the particular value of $\theta$ used in (2.3). We also use the subscript $\theta$ with certain other
quantities to emphasize their dependence on $\theta.$

We recall [7,11,14-15] that the bound states for the Schr\"odinger operator
associated with (2.1)-(2.3) correspond to
square-integrable solutions to (2.1) satisfying the boundary condition
(2.3). Therefore, the bound-state energies, i.e. the $k^2$-values at which bound states occur, depend on the boundary parameter $\theta.$
When $V$ belongs to class $\Cal A$ given in (2.2),
it is known [7,11,14-16] that there can be at most a finite number of
bound states and that the number of bound states is also
affected by the parameter $\theta.$
Because of the
selfadjointness of the corresponding Schr\"odinger operator,
each bound-state energy must be real. It is already
known [7,11,14-16] that for each positive $k^2$-value in (2.1) there
correspond two linearly independent solutions, e.g. $f(k,x)$
and $f(-k,x),$ neither of which is square integrable in $x\in\bR^+$
as a result of (2.4). Each
bound state is known [7,11,14-16] to be simple in the sense that there
exists only one linearly independent square-integrable solution
to (2.1) satisfying (2.3) at a bound-state energy.
The bound states, if there are any, can only
occur at certain negative values of $k^2,$ and we will assume
that they occur at $k=i\gamma_s$ for $s=1,\dots,N$ for some
nonnegative integer $N$ and distinct positive values $\gamma_s.$
Note that the $\gamma_s$-values are not
in an increasing or decreasing order.
%The order of $\gamma_s$ is such that
%we first add the bound state at $k=i\gamma_1,$ then add the bound state
%at $k=i\gamma_2,$ and so on. We remove the bound states in the reverse order,
%i.e. we first remove the bound state at $k=i\gamma_N,$ then remove
%the bound state at $k=i\gamma_{N-1},$ and so on.
Note also that even though the
value of $N$ and the values of $\gamma_s$ all depend
on the choice of $\theta,$ for notational simplicity we usually suppress the
dependence on $\theta$ for those quantities.

The so-called Jost function associated with
(2.1)-(2.3), usually denoted by $F_\theta(k),$ is defined [7,11,14,15] as
$$F_\theta(k):=\cases -i[f'(k,0)+\cot\theta\,f(k,0)],&\qquad \theta\in (0,\pi),\\
\stretch
f(k,0),&\qquad \theta=\pi,
\endcases\tag 2.6$$
and it helps us to identify the bound states and to define
the scattering matrix.
It is known [8-10,14-16] that $f(k,x)$ and $f(-k,x)$
are linearly independent for each fixed $k\in\bC\setminus\{0\}.$ Thus,
we can express the regular solution $\varphi_\theta(k,x)$
appearing in (2.5) as a linear combination of
$f(k,x)$ and $f(-k,x).$ In fact, with the help of (2.5) and
(2.6) we get
$$\varphi_\theta(k,x)=\cases
\ds\frac{1}{2k}\left[
F_\theta(k)\,f(-k,x)-F_\theta(-k)\,f(k,x)\right],&\qquad \theta\in
(0,\pi),\\
\stretch
\ds\frac{i}{2k}\left[
F_\theta(k)\,f(-k,x)-F_\theta(-k)\,f(k,x)\right],&\qquad \theta=
\pi.\endcases\tag 2.7$$
 From (2.3), (2.4), and (2.7) we see that a bound state can only occur
at a zero of $F_\theta(k),$ which is equivalent to the linear dependence of
the two solutions $\varphi_\theta(k,x)$ and $f(k,x)$ at that particular $k$-value.
This is because the linear dependence on
$\varphi_\theta(k,x)$ assures the satisfaction of the boundary condition (2.3), and
the linear dependence on $f(k,x)$ guarantees an exponential decay
as $x\to+\infty$ and in turn the square integrability in $x\in\bR^+.$

We have seen that there are at most a finite number of zeros of the
Jost function $F_\theta(k)$ in $\bCp$ and such zeros can only occur
on the positive imaginary axis, and those zeros correspond to bound states of the
Schr\"odinger operator given in (2.1)-(2.3).
Let us now consider the zeros of $F_\theta(k)$ in $\bCm,$
which are called resonances. When
$V(x)\equiv 0,$ from (2.4) and (2.6) it follows that
$$F_\theta(k)=\cases k-i\,\cot\theta,\qquad &\theta\in(0,\pi),\\
\stretch
1,\qquad & \theta=\pi.\endcases\tag 2.8$$
Thus, the number of resonances is at most one when
$V\equiv 0.$ As stated in Theorem~2.1(g) later,
if $V\not\equiv 0$ then there must be a countably infinite
number of resonances, and each resonance
occurs either on the negative imaginary axis or a pair of resonances are symmetrically
located with respect to the negative imaginary axis.

In our paper we are primarily
interested in imaginary resonances, i.e. those resonances located on the negative
imaginary axis. Through a pathological example [19] it is known that the number of
imaginary resonances can be countably infinite even when the potential $V$ is in class $\Cal A.$
On the other hand, the number of imaginary resonances is guaranteed to be
finite under some minimal further assumptions, e.g.
see Proposition~7 of [19], such as $V(x)\ge 0$ or $V(x)\le 0$
in some neighborhood of $x=b,$ where $b$ is the parameter appearing in (2.2)
and related to the compact support
of $V.$
In Section~3 we develop various equivalent criteria to identify
each imaginary resonance either as an eligible resonance or an ineligible
resonance and explore the connection between bound states
and eligible resonances.

Having seen that the zeros of $F_\theta(k)$ in $\bCp$ correspond to bound states and
the zeros in $\bCm$ correspond to resonances, let us now consider zeros of
$F_\theta(k)$ occurring on the real axis. It is known [7,14,15] that the only real zero of $F_\theta(k)$ can occur at $k=0$
and such a zero, if it exists, must be a simple zero.
The case $F_\theta(0)=0$
corresponds to the exceptional case, and the case $F_\theta(0)\ne 0$ corresponds to
the generic case. In the exceptional case, the number of bound states
may change by one under a small perturbation of the potential.
Let us also consider the Jost solution $f(k,x)$ and
the regular solution $\varphi_\theta(k,x)$ appearing in (2.4) and (2.5),
respectively, at $k=0.$
Generically $\varphi_\theta(0,x)$ becomes
unbounded as $x\to+\infty,$ whereas in the exceptional case
it remains bounded as $x\to+\infty.$ The behavior of
$\varphi_\theta(0,x)$ as $x\to+\infty$ is obtained by letting $k\to 0$ in (2.7),
using (2.4), and exploiting the known
behaviors of $f(0,x)$ and $\dot f(0,x)$ as
$x\to+\infty,$ where we use an overdot to indicate the $k$-derivative. As
seen from (2.4) we have $f(0,x)\equiv 1$ and $\dot f(0,x)=ix$
for $x\ge b.$ From (2.7) at $k=0$ we get
$$\varphi_\theta(0,x)=\cases \dot F_\theta(0)\,f(0,x)-F_\theta(0)\,\dot f(0,x),
\qquad & \theta\in (0,\pi),\\
\stretch
i\left[\dot F_\theta(0)\,f(0,x)-F_\theta(0)\,\dot f(0,x)
\right],
\qquad & \theta=\pi),\endcases$$
which shows that
$\varphi_\theta(0,x)$ is proportional to
$f(0,x)$ and hence remains bounded in the exceptional case and
that $\varphi_\theta(0,x)$ contains $\dot f(0,x)$ and hence becomes unbounded
in the generic case.

Recall that we assume the bound states occur at the zeros
$k=i\gamma_s$ of $F_\theta(k)$ appearing
in (2.6) for $s=1,\dots,N.$ It is known [6,14,15] that $\varphi_\theta(i\gamma_s,x)$ is
real valued and square integrable. The positive
quantity $g_s$ defined as
$$g_s:=\ds\frac{1}{\ds\sqrt{\int_0^\infty dx\,\varphi_\theta(i\gamma_s,x)^2}},
\qquad s=1,\dots,N,\tag 2.9$$
is known as the Gel'fand-Levitan norming constant
for the bound state at $k=i\gamma_s.$
Let us use
$\Cal G$ to denote the Gel'fand-Levitan spectral data set [7,8,14,15] given by
$$\Cal G:=\big\{|F_\theta(k)|:\ k\in\bR;\ \{\gamma_s,g_s\}_{s=1}^N\big\}.\tag 2.10$$
We refer to the information consisting of $|F_\theta(k)|$
for $k\in\bR$ as the continuous part of the Gel'fand-Levitan spectral data and
refer to the portion $\{\gamma_s,g_s\}_{s=1}^N$ as the discrete part of
the Gel'fand-Levitan
spectral data. For the construction of $V$ and $\theta$ from $\Cal
G$ via the Gel'fand-Levitan method, we outline
the recovery procedure below and refer the reader to
[7,11,14,15] for the details.

\item{(a)} From the large-$k$ asymptotics [7]
$$|F_\theta(k)|=\cases |k|+O(1),&\qquad  k\to\pm\infty,\quad
\theta\in(0,\pi),\\
\stretch 1+O\left(\ds\frac{1}{k}\right),&\qquad k\to\pm\infty,\quad
\theta=\pi,\endcases\tag 2.11$$
we can tell whether
$\theta\in(0,\pi)$ or $\theta=\pi.$

\item{(b)} We form [7,11,14,15] the Gel'fand-Levitan kernel
$G_\theta(x,y),$
where for $\theta\in(0,\pi)$ we have
$$G_\theta(x,y):=\ds\frac{1}{\pi}\int_{-\infty}^\infty
dk\,\left[\ds\frac{k^2}{|F_\theta(k)|^2}-1 \right](\cos kx)(\cos
ky)+\ds\sum_{s=1}^N g_s^2\,(\cosh\gamma_s x)(\cosh\gamma_s y),\tag 2.12$$
 and for $\theta=\pi$ we have
$$G_\theta(x,y):=
\ds\frac{1}{\pi}\int_{-\infty}^\infty
dk\,\left[\ds\frac{1}{|F_\theta(k)|^2}-1
\right](\sin kx)(\sin ky)+\ds\sum_{s=1}^N \ds\frac{g_s^2}
{\gamma_s^2}\,(\sinh\gamma_s x)(\sinh\gamma_s y).\tag 2.13$$

\item{(c)} Using $G_\theta(x,y)$ as input to the
Gel'fand-Levitan integral equation
$$A_\theta(x,y)+G_\theta(x,y)+\int_0^x
dz\,A_\theta(x,z)\,G_\theta(z,y),\qquad 0<y<x,\tag 2.14$$
we
obtain $A_\theta(x,y).$ The unique solvability of (2.14) is
known [11,14,15] for the spectral data set corresponding to
a potential in class $\Cal A$
and a boundary condition as in (2.3).

\item{(d)} We obtain the potential $V(x)$ and the
boundary parameter $\theta$ via [7,11,14,15]
$$V(x)=2\ds\frac{d}{dx}A_\theta(x,x),\qquad \theta\in(0,\pi],\tag 2.14$$
$$\cot\theta=-A_\theta(0,0),\qquad \theta\in(0,\pi).$$

\item{(e)} The regular solution $\varphi_\theta(k,x)$
is recovered from $A_\theta(x,y)$ via [7,11,14,15]
$$\varphi_\theta(k,x)=\cases \cos kx+\ds\int_0^x dy\, A_\theta(x,y)\,
\cos ky,\qquad & \theta\in(0,\pi),\\
\stretch
\ds\frac{\sin kx}{k}+\ds\int_0^x dy\, A_\theta(x,y)\,
\ds\frac{\sin ky}{k},\qquad & \theta=\pi.\endcases$$

An alternative to the Gel'fand-Levitan procedure is the Marchenko
method [7,14,15], which uses the input data set
$\Cal M$ given by
$$\Cal M:=\{S_\theta(k):\ k\in\bR;\ \{\gamma_s,m_s\}_{s=1}^N\},\tag 2.15$$
where the scattering matrix $S_\theta(k)$ is defined in
terms of the Jost function $F_\theta(k)$ as [7,14,15]
$$S_\theta(k):=\cases -\ds\frac{F_\theta(-k)}{F_\theta(k)},&\qquad \theta\in (0,\pi),\\
\stretch
\ds\frac{F_\theta(-k)}{F_\theta(k)},&\qquad \theta=\pi,\endcases\tag 2.16$$
and the
Marchenko bound-state norming constants $m_s$ are given by [7,14,15]
$$m_s:=\ds\frac{1}{\ds\sqrt{\int_0^\infty dx\,f(i\gamma_s,x)^2}},
\qquad s=1,\dots,N.\tag 2.17$$
We refer to the information consisting of $S_\theta(k)$ for $k\in\bR$ as the continuous
part of the Marchenko scattering data
and the portion $\{\gamma_s,m_s\}_{s=1}^N$ as the
discrete part of the scattering data.

For the construction of $V$ and $\theta$ from $\Cal M$ given in (2.15), we
outline the steps of the
Marchenko recovery method below and
refer the reader to [7,14,15] for further details.

\item{(a)} Using the data $\Cal M,$ we construct the Marchenko
kernel $M_\theta$ as
$$M_\theta(y):=\cases
\ds\frac{1}{2\pi} \int_{-\infty}^\infty dk\, \left[S_\theta(k)-1\right] e^{iky}+
\sum_{s=1}^N m_s^2\, e^{-\gamma_s y},&\qquad \theta\in(0,\pi),\\
\ds\frac{1}{2\pi} \int_{-\infty}^\infty dk\, \left[1-S_\theta(k)\right] e^{iky}+
\sum_{s=1}^N m_s^2\, e^{-\gamma_s y},&\qquad \theta=\pi.
\endcases\tag 2.18$$

\item{(b)} Using $M_\theta(y)$ given in (2.18) as input to the Marchenko
integral equation
$$K(x,y)+M_\theta(x+y)+\ds\int_x^\infty dz\,K(x,z)\,M_\theta(z+y)=0,\qquad y>x,\tag 2.19$$
we obtain $K(x,y).$ The unique solvability of (2.19) is guaranteed [8-10,14,15]
if the scattering data set corresponds to
a potential in class $\Cal A$ given in (2.2).

\item{(c)} The potential $V(x)$ and the Jost solution $f(k,x)$ are
obtained from $K(x,y)$ via
$$V(x)=-2\,\ds\frac{dK(x,x)}{dx},\quad f(k,x)=e^{ikx}+\int_x^\infty
dy\,K(x,y)\,e^{iky}.\tag 2.20$$

\item{(d)} Having $K(x,y)$ and $S_\theta(k)$ at hand, we can
recover $\cot\theta$ as well. For this purpose, we can
proceed as follows. From the second equation in (2.20) we get
$$f(k,0)=1+\int_0^\infty dy\,K(0,y)\,e^{iky},\tag 2.21$$
$$f'(k,0)=ik-K(0,0)+\int_0^\infty dy\,K_x(0,y)\,e^{iky},\tag 2.22$$
where $K_x(0,y)$ denotes the $x$-derivative of $K(x,y)$ evaluated at $x=0.$
In light of the second line of (2.16) we then check if we have
$$S_\theta(k)=\ds\frac{1+\ds\int_0^\infty dy\,K(0,y)\,e^{-iky}}
{1+\ds\int_0^\infty dy\,K(0,y)\,e^{iky}},
\tag 2.23$$
which is obtained by using (2.21) and (2.22) in the second line of (2.16).
We conclude that $\theta=\pi$ if (2.23) is satisfied.
If (2.23) is not satisfied, we conclude that $\theta\in(0,\pi)$
and uniquely determine $\cot\theta$ as
$$\cot\theta=\ds\frac{-f'(-k,0)-S_\theta(k)\,f'(k,0)}{f(-k,0)+S_\theta(k)\,f(k,0)},\tag 2.24$$
which is obtained with the help of (2.6), (2.16), (2.21), and (2.22).

For easy citation later on, we summarize the results presented above and several additional known facts [1,7-10,14-16] in the following theorem.

\noindent {\bf Theorem 2.1} {\it Consider the Schr\"odinger
operator given in (2.1)-(2.3) with the potential $V$
in class $\Cal A,$ a fixed boundary parameter $\theta\in(0,\pi],$
and $b$ being the constant appearing in (2.2)
related to the compact support of the potential.
Let $F_\theta(k)$ be the corresponding Jost function given in (2.6) and
$S_\theta(k)$ be the corresponding scattering matrix appearing in (2.16).
Then:}

\item{(a)} {\it The Jost function $F_\theta(k)$ has an analytic
extension from $k\in\bR$ to the entire complex plane $\bC.$ There are at most a
finite number of zeros of $F_\theta(k)$ in $\bCp,$ they
occur on the positive imaginary axis, say at
$k=i\gamma_s$ for $s=1,\dots,N,$ they are all simple,
and they correspond to the bound states of (2.1) with the selfadjoint
boundary condition (2.3).
A real zero of $F_\theta(k)$
can only occur at $k=0,$ and such a zero, if it exists, must be simple.}

\item{(b)} {\it As $k\to\infty$ in $\bCpb$ we have}
$$F_\theta(k)=\cases
k-i\cot\theta+\ds\frac{i}{2}\int_0^b dx\,V(x)+o(1),&\qquad \theta\in(0,\pi),\\
\stretch
1-\ds\frac{1}{2ik}\int_0^b dx\,V(x)+o\left(\ds\frac{1}{k}\right),&\qquad \theta=\pi.\endcases
\tag 2.25$$

\item{(c)} {\it As $k\to\infty$ in $\bCmb$ we have}
$$F_\theta(k)=\cases
k-i\cot\theta+\ds\frac{i}{2}\int_0^b dx\,V(x)
+e^{2ikb}o(1),&\qquad \theta\in(0,\pi),\\
\stretch
1-\ds\frac{1}{2ik}\int_0^b dx\,V(x)+e^{2ikb}o\left(\ds\frac{1}{k}\right),&\qquad \theta=\pi.
\endcases$$

\item{(d)} {\it As $k\to\pm\infty$ in ${\bold R},$ the large-$|k|$ asymptotics of the scattering matrix
$S_\theta(k)$ is given by}
$$S_\theta(k)=\cases
1-\displaystyle\frac{i}{k}\int_0^b dx\,V(x)
+\displaystyle\frac{2i}{k}\,\cot\theta+o\left(\displaystyle\frac{1}{k}\right),&\qquad \theta\in(0,\pi)
,\\
\stretch
1-\displaystyle\frac{i}{k}\int_0^b dx\,V(x)+o\left(\displaystyle\frac{1}{k}\right)
,&\qquad \theta=\pi
.
\endcases$$

\item{(e)} {\it The scattering matrix $S_\theta(k)$ defined in (2.9)
has a meromorphic extension from $k\in\bR$
to $k\in\bC.$
The poles of $S_\theta(k)$ in $\bCp$ are all simple
and occur at
$k=i\gamma_s$ for $s=1,\dots,N.$ The Marchenko
norming constants
$m_s$ defined in (2.17) are related to the residues of the
scattering matrix at those poles as}
$$\text{Res}(S_\theta,i\gamma_s)=\cases i\,m_s^2,
&\qquad \theta\in(0,\pi),\\
\stretch
-i\,m_s^2,
&\qquad \theta=\pi,\endcases\tag 2.26$$
{\it where
$\text{\rm Res}(S_\theta,i\gamma_s)$ denotes the residue of
$S_\theta(k)$ at $k=i\gamma_s.$}

\item{(f)} {\it For each $\theta\in (0,\pi],$ the scattering
matrix $S_\theta(k)$ is analytic at $k=0$ in $\bC.$
The value of $S_\theta(0)$ is either $+1$ or $-1.$ Specifically,
for $\theta=\pi$ we have}
$$S_\pi(0)=\cases +1, &\qquad f(0,0)\ne 0,\\
\stretch
-1, &\qquad f(0,0)= 0
,\endcases\tag 2.27$$
{\it and for any $\theta\in(0,\pi)$ we have}
$$S_\theta(0)=\cases -1, &\qquad F_\theta(0)\ne 0,\\
\stretch
+1, &\qquad F_\theta(0)= 0
.\endcases\tag 2.28$$

\item{(g)} {\it Unless $V(x)\equiv 0,$ there are infinitely
many zeros of $F_\theta(k)$ in ${\bold C^-},$ and such zeros are known as resonances.
The resonances need not be simple,
and they are located either on the negative imaginary axis
or occur in pairs located symmetrically with respect to the negative imaginary axis.}

\item{(h)} {\it The Gel'fand-Levitan norming constants $g_s$ appearing in
(2.9) and the Marchenko norming constants $m_s$ appearing in (2.17) are related to each other
as}
$$g_s=\ds\frac{2\gamma_s \, m_s}{|F_{\theta}(-i\gamma_s)|}, \qquad
\theta\in(0,\pi].
\tag 2.29$$

\item{(i)} {\it The potential $V$ and the boundary parameter
$\theta$ are uniquely determined from the Gel'fand-Levitan spectral data
$\Cal G$ given in (2.10).}

\item{(j)} {\it The potential $V$ and the boundary parameter
$\theta$ are uniquely determined from the Marchenko scattering data
$\Cal M$ given in (2.15).}

\noindent PROOF: For (a), (i), (j), we refer the reader to
[7,14,15]. For (b), (c), (d), (e), the reader is referred to [6]. The result in (f)
is obtained by using (2.6) and (2.16)
with the help of a series expansion around $k=0.$
The proof of (g) is as follows. From (2.8) we already
know that the number of resonances corresponding to $V(x)\equiv 0$
is either zero or one. For $V(x)\not\equiv 0$ with $b>0$ in (2.2),
we conclude, from (a)-(c), that $e^{2ikb} F_\theta(k)$
is entire in $k$ and behaves as
$O(k)$ as $k\to\infty$ in $\bC.$ If $F_\theta(k)$
had no zeros or had only a finite number of zeros in $\bC,$
then the Hadamard
factorization of $e^{2ikb} F_\theta(k)$ and the use of Liouville's theorem would force $F_\theta(k)$
to be equal to $e^{-2ikb}$ multiplied with either
a constant
or a polynomial in $k.$ However,
such a behavior would contradict (2.25).
Thus, the number of resonances must be countably infinite.
Since $k$ appears as $ik$ in
$f(k,0)$ and $f'(k,0),$ it follows from (2.6), that the zeros
of $F_\theta(k)$ in $\bCm$ either occur on the negative
imaginary axis or a pair of resonances are symmetrically located with respect to the negative imaginary axis. From Example~6.2(c) we know that a resonance need not be simple. Thus, the proof of (g) is complete.
 Note that (2.29) holds for $\theta=\pi$ as well for $\theta\in(0,\pi).$
The result in (2.29) is obtained by evaluating (2.7) at
the bound state $k=i\gamma_s,$ using $F_\theta(i\gamma_s)=0$ in that equation,
taking the square of both sides of the resulting equation,
followed by an integration on $x\in\bR^+,$ and finally by using
(2.9) and (2.17) in the resulting equation. \qed

Next, we elaborate on the exceptional case for the half-line Schr\"odinger
operator and present the behavior of the corresponding
scattering coefficients for the full-line
Schr\"odinger
operator at $k=0.$
Such results are needed in Sections~3 and 4 in the elaboration of the
nonuniqueness arising in the special case, i.e. case (iii) of Section~4.

Recall that the exceptional case for the half-line Schr\"odinger
operator occurs when $F_\theta(0)=0,$ where $F_\theta(k)$
is the Jost function defined in (2.6). Since we can view the potential $V$
appearing in (2.1) as the potential on the full line with
$V(x)\equiv 0$ for $x<0,$ we can uniquely [7,8] associate with
$V$ the scattering coefficients $T,$ $L,$ $R,$ where
$T$ is the transmission coefficient, $L$ is the reflection coefficient from the
left, and $R$ is the reflection coefficient from the
right. This is done via [7,8]
$$f(k,0)=\ds\frac{1+L(k)}{T(k)},\quad
f'(k,0)=ik\, \ds\frac{1-L(k)}{T(k)},\quad
R(k)=-\ds\frac{L(-k)\,T(k)}{T(-k)}.\tag 2.30$$
The exceptional case for the full-line Schr\"odinger
operator occurs when $T(0)\ne 0,$ and the generic case
occurs when $T(0)=0.$

\noindent {\bf Theorem 2.2} {\it Consider the half-line Schr\"odinger
operator given in (2.1)-(2.3) with the potential $V$
in class $\Cal A$ and with a fixed boundary parameter $\theta\in(0,\pi].$
Let $f(k,x)$ and
$F_\theta(k)$ be the corresponding Jost solution and the Jost
function appearing in (2.4) and (2.6), respectively.
Further, let $T(k),$ $L(k),$ $R(k)$ be the corresponding scattering coefficients
appearing in (2.30).
Then:}

\item{(a)} {\it The half-line exceptional case with the Dirichlet boundary
condition, i.e. $f(0,0)=0,$ corresponds to the following zero-energy
behavior of the scattering coefficients:}
$$T(0)=0,\quad \dot T(0)\ne 0,\quad  L(0)=-1,\quad \dot L(0)=0
,\quad \ddot L(0)\ne 0,\tag 2.31$$
$$R(0)=-1,\quad \dot R(0)=-\ds\frac{\ddot T(0)}{\dot T(0)}
,\quad \ddot R(0)=-\dot T(0)^2-\ds\frac{\ddot T(0)^2}{\dot T(0)^2} ,\tag 2.32$$
{\it where we recall that an overdot denotes the $k$-derivative.}

\item{(b)} {\it The half-line exceptional case with the Neumann boundary
condition, i.e. $f'(0,0)=0,$ corresponds to the following zero-energy
behavior of the scattering coefficients:}
$$T(0)\ne 0,\quad L(0)\ne -1
,\quad R(0)\ne -1.\tag 2.33$$

\item{(c)} {\it The half-line exceptional case with the non-Dirichlet and non-Neumann boundary
conditions, i.e. $F_\theta(0)=0$ with
$\theta\in(0,\pi/2)\cup (\pi/2,\pi),$
corresponds to the following zero-energy
behavior of the scattering coefficients:}
$$T(0)=0,\quad \dot T(0)\ne 0,\quad  L(0)=-1,\quad R(0)=-1
,$$
$$\dot L(0)=-\ds\frac{2i}{\cot\theta}
%,\quad \ddot L(0)=\ds\frac{4}{\cot^2\theta}
,\quad \dot R(0)=\ds\frac{2i}{\cot\theta}-\ds\frac{\ddot T(0)}{\dot T(0)}.\tag 2.34$$

\noindent PROOF: The behavior of the scattering coefficients around $k=0$ is
already known [5,8]. In the full-line generic case we have
$$T(0)=0,\quad \dot T(0)\ne 0,\quad L(0)=-1,\quad R(0)=-1,\tag 2.35$$
and in the full-line exceptional case we have
$$T(0)\ne 0,\quad L(0)\in(-1,1),\quad R(0)\in(-1,1).\tag 2.36$$
 From Theorem~2.1(a), when $V\in\Cal A$ we know that $f(k,0)$ and $f'(k,0)$ are entire, and hence
with the help of (2.30) we see that $T(k),$ $R(k),$ and $L(k)$ are analytic at $k=0.$
Expanding around $k=0$
the first identity in (2.30), we see that (2.36) is incompatible
with $f(0,0)=0$ and hence in case of (a) in our theorem, we must have (2.35). Then,
the expansion of the first identity in (2.30) yields
$$f(0,0)+k\,\dot f(0,0)+O(k^2)=\ds\frac{\dot L(0)}{\dot T(0)}+\ds\frac{k}{2}\left[
\ds\frac{\ddot L(0)}{\dot T(0)}-\ds\frac{\dot L(0)\,\ddot T(0)}{\dot T(0)^2}\right]
+O(k^2),
\qquad k\to 0
\text{ in } \bC.\tag 2.37$$
 From Theorem~2.1(a) we already know that $k=0$ must be a simple zero of
 $f(k,0)$ and hence $\dot f(0,0)\ne 0.$ Thus, from (2.37) we get
 $\dot L(0)=0$ and $\ddot L(0)\ne 0.$ Hence, we have proved
 (2.31). In fact, the expansion around $k=0$ of
 the identity [8-10]
 $$L(k)\,L(-k)+T(k)\,T(-k)=1,\qquad k\in\bC,$$
indicates that in the full-line generic case we have
$$\ddot L(0)+\dot L(0)^2+\dot T(0)^2=0,\tag 2.38$$
and hence (2.38) shows that in case of (a) we have
$$\ddot L(0)=-\dot T(0)^2,\tag 2.39$$
which also confirms that $\ddot L(0)\ne 0$ in (2.31).
We establish (2.32), by expanding around $k=0$ the third identity in (2.30)
and using (2.31) and (2.39). Let us now turn to the proof of (b).
Expanding around $k=0$ the second identity in (2.30), we see that
(2.35) is incompatible with $f'(0,0)=0.$ Thus, we must have (2.36) in case of (b),
which establishes (2.33). Finally, let us prove (c). Using the first two
identities in
(2.6), we get
$$F_\theta(k)=k\,\ds\frac{1-L(k)}{T(k)}-i\,\cot\theta\,\ds\frac{1+L(k)}{T(k)}.
\tag 2.40$$
Note that (2.36) is not compatible with $\cot\theta\ne 0$ and $F_\theta(0)=0.$
Thus, we must have (2.35) in case of (c). Then, expanding around $k=0$ both sides of
(2.40) we get
$$F_\theta(0)=\ds\frac{2-i\,\cot\theta\, \dot L(0)}{\dot T(0)},\tag 2.41$$
$$\dot F_\theta(0)=-
\ds\frac{\ddot T(0)}{2\,\dot T(0)}\,\ds\frac{2-i\,\cot\theta\, \dot L(0)}{\dot T(0)}-\ds\frac{2 \,\dot L(0)+i\,\cot\theta\,\ddot L(0)}{2\,\dot T(0)}
.$$
Since $F_\theta(0)=0,$ from (2.41) we get $\dot L(0)=-2i/\cot\theta.$ Finally,
with the help of (2.30) we get $\dot R(0)$ given in (2.34).
\qed

The next theorem shows that
if the half-line
Schr\"odinger
operator with the Neumann boundary condition
and with a potential $V$ belonging to class
$\Cal A$ has no bound states then the full-line Schr\"odinger
operator with the same potential $V$ cannot have any bound states either. The result
is needed for the proof of Theorem~2.4 and later in the analysis
in Section~4.

\noindent {\bf Theorem 2.3} {\it Consider the half-line Schr\"odinger
operator given in (2.1)-(2.3) with the potential $V$
in class $\Cal A$ and with a fixed boundary parameter $\theta\in(0,\pi],$
and let $f(k,x)$ and
$F_\theta(k)$ be the corresponding Jost solution and the Jost
function appearing in (2.4) and (2.6), respectively. Let $N_\theta$ denote the number of bound states, i.e. the number of zeros of
$F_\theta(i\beta)$ when $\beta\in(0,+\infty).$
Let $T(k),$ $L(k),$ $R(k)$ be the corresponding scattering coefficients
appearing in (2.30). Let $\tilde N$ denote the number of bound states for the
corresponding
full-line Schr\"odinger operator, i.e. let $\tilde N$ denote the number of zeros of $1/T(i\beta)$ in the interval $\beta\in(0,+\infty).$
If $N_{\pi/2}=0$ then we must have $\tilde N=0.$}

\noindent PROOF: It is already known [7]
that $N_{\theta_1}\le N_{\theta_2}$ if $\theta_1\ge \theta_2.$
Thus, in particular we have
$N_\pi\le N_{\pi/2}.$ Since we assume $N_{\pi/2}=0,$ we then
also have $N_\pi=0.$ Thus, neither $f'(i\beta,0)$ nor $f(i\beta,0)$ vanishes
for $\beta>0.$ From (2.25) we then conclude that
$-f'(i\beta,0)>0$ and $f(i\beta,0)>0$ for all $\beta>0.$
The first two identities in (2.30) yield
 $$\ds\frac{2ik}{T(k)}=f'(k,0)+ik\, f(k,0),\qquad k\in\bC.\tag 2.42$$
 From (2.42), using $k=i\beta$ we obtain
$$\ds\frac{2\beta}{T(i\beta)}=-f'(i\beta,0)+\beta\,f(i\beta,0).\tag 2.43$$
Since the right-hand side of (2.43) is positive for all $\beta>0,$ we conclude
that $T(i\beta)$ does not have any poles for $\beta>0$ and hence
$\tilde N=0.$ \qed

The following theorem shows that in the absence of any bound states, the Marchenko
equation given in (2.19) is equivalent to the full-line
Marchenko equation given by
$$K(x,y)+\hat R(x+y)+\ds\int_x^\infty dz\,K(x,z)\,\hat R(z+y)=0,\qquad y>x,\tag 2.44$$
where $\hat R(y)$ denotes the Fourier transform of the reflection coefficient
$R(k)$ appearing in (2.30), namely
$$\hat R(y):=\ds\frac{1}{2\pi}\ds\int_{-\infty}^\infty dk\,R(k)\,e^{iky}.\tag 2.45$$
The result in Theorem~2.4 is needed in the characterization of the
double nonuniqueness in the special case in Section~4, i.e. case (iii) there.

\noindent {\bf Theorem 2.4} {\it Consider the half-line Schr\"odinger
operator given in (2.1)-(2.3) with the potential $V$
in class $\Cal A$ and with a fixed boundary parameter $\theta\in(0,\pi].$
Let $F_\theta(k),$ $S_\theta(k),$ and $M_\theta(y)$
 be the corresponding Jost function, the
scattering matrix, and the Marchenko kernel defined
in (2.6), (2.16), and (2.18), respectively.
Let $T(k),$ $L(k),$ $R(k)$ be the corresponding scattering coefficients
appearing in (2.30). Assume that neither the half-line Schr\"odinger
operator nor the full-line Schr\"odinger
operator has any bound states,
i.e. $F_\theta(k)$ has no zeros
on the positive imaginary axis and $T(k)$ has no poles
on the positive imaginary axis.
Then, we have}
$$M_\theta(y)=\hat R(y), \qquad y>0,\quad \theta\in(0,\pi],
\tag 2.46$$
{\it where $\hat R(y)$ is the quantity given in (2.45).}

\noindent PROOF: From (2.30) we get
$$\ds\frac{2ik\,L(k)}{T(k)}=ik\,f(k,0)-f'(k,0),\tag 2.47$$
and hence from (2.42) and (2.47) we have
$$T(k)=\ds\frac{2ik}{f'(k,0)+ik\,f(k,0)},\qquad
\ds\frac{L(-k)}{T(-k)}=\ds\frac{f'(-k,0)+ik\,f(-k,0)}{2ik}.\tag 2.48$$
Using (2.48) in the third equation in (2.30) we obtain
$$R(k)=-\ds\frac{f'(-k,0)+ik\,f(-k,0)}{f'(k,0)+ik\,f(k,0)}.\tag 2.49$$
Thus, from (2.49) and the second line of (2.16) we get
$$1-S_\pi(k)-R(k)=1-\ds\frac{f(-k,0)}{f(k,0)}+\ds\frac{f'(-k,0)
+ik\,f(-k,0)}{f'(k,0)+ik\,f(k,0)}.\tag 2.50$$
Using the Wronskian relation [8-10]
$$f(-k,0)\,f'(k,0)-f'(-k,0)\,f(k,0)=2ik,$$
and the first equality in (2.48), we can rewrite (2.50) as
$$1-S_\pi(k)-R(k)=1-\ds\frac{T(k)}{f(k,0)}.\tag 2.51$$
In the absence of bound states for the full-line
Schr\"odinger equation, it is known [8-10] that $T(k)$
is analytic in $\bCp$ and continuous in $\bCpb$ and
$T(k)=1+O(1/k)$ as $k\to\infty$ in $\bCpb.$
In the absence of bound states for the half-line
Schr\"odinger equation, $f(k,0)$ and hence also $1/f(k,0)$ are
analytic in $\bCp$ and continuous in $\bCpb$ and
behave as $1+O(1/k)$ as $k\to\infty$ in $\bCpb.$
Furthermore, from Theorem~2.2(a) the continuity of $T(k)/f(k,0)$
at $k=0$ is assured.
Thus, the right-hand side of (2.51) is analytic
in $\bCp$ and continuous in $\bCpb$ and behaves
as $O(1/k)$ as $k\to\infty$ in $\bCpb.$ Hence, its Fourier transform vanishes
for $y>0,$ i.e.
$$\ds\frac{1}{2\pi}\int_{-\infty}^\infty dk\,\left[1-S_\pi(k)-R(k)\right]\,
e^{iky}=0,\qquad y>0.\tag 2.52$$
Comparing (2.52) with (2.45) and the second line of (2.18) without the summation
term there, we see that $M_\pi(y)=\hat R(y)$ for $y>0,$ establishing (2.46)
for $\theta=\pi.$
In a similar way, we can show that, for $\theta\in(0,\pi),$ we have
$$R(k)-S_\theta(k)+1=1-\ds\frac{(k+i\,\cot\theta)\,T(k)}{F_\theta(k)}.\tag 2.53$$
In the absence of bound states for the half-line Schr\"odinger operator, from
Theorem~2.1 we know that $1/F_\theta(k)$ is analytic in $\bCp,$ continuous
in $\bCpb\setminus\{ 0\},$ and behaves like $O(1/k)$ as
$k\to\infty$ in $\bCpb.$ In the absence of bound states for the full-line
Schr\"odinger operator, we already know
that $T(k)$
is analytic in $\bCp,$ continuous in $\bCpb,$ and behaves as
$1+O(1/k)$ as $k\to\infty$ in $\bCpb.$
Furthermore, from (b) and (c) of Theorem~2.2 it follows that the
second term on the right-hand side in (2.53)
is continuous at $k=0.$
Thus, the right-hand side in (2.53) is analytic
in $k\in\bCp$ and continuous in $k\in\bCpb$ and behaves as $O(1/k)$ as
$k\to\infty$ in $\bCpb.$ Hence, its Fourier transform for $y>0$ vanishes, i.e.
we have
$$\ds\frac{1}{2\pi}\int_{-\infty}^\infty dk\,\left[R(k)-S_\theta(k)+1\right]\,
e^{iky}=0,\qquad y>0,\tag 2.62$$
yielding $M_\theta(y)=\hat R(y)$ for $y>0.$ Therefore, (2.46) holds also when
$\theta\in(0,\pi).$ \qed

\vskip 10 pt
\noindent {\bf 3. DARBOUX TRANSFORMATION AND ELIGIBLE RESONANCES}
\vskip 3 pt

Recall that a Darboux transformation [8,9,15] allows us to change
the discrete spectrum of a differential operator by adding or removing
a finite number of discrete eigenvalues without changing
the continuous spectrum. In preparation for the analysis in Section~5, in this section we
provide the Darboux transformation formulas when a bound state is added
or removed from the spectrum of the Schr\"odinger operator on the
half line. We also provide various results related to the Darboux
transformation with
compactly-supported potentials. In particular, we provide
the necessary and sufficient conditions for retaining the compact-support property
of the potential when we add a bound state. We show that such a bound state
can only come from an eligible resonance,
which is a zero of the Jost function $F_\theta(k)$ occurring on the
negative imaginary axis and can be converted to a bound state via a Darboux
transformation without changing the compact support of the potential
satisfying a certain derivative condition.
We provide various equivalent characterizations of eligible resonances,
such as (3.19), (3.38), and (3.53).

For clarity, we use the notation
$\theta_j,$ $V(x;j),$ $\varphi(k,x;j),$ and
$F(k;j)$ to denote the relevant
quantities corresponding to the Schr\"odinger operator
with bound states at $k=i\gamma_1,\dots,i\gamma_j,$
where
the case $j=0$ refers to the quantities
without bound states. Note that $\theta_j$ is the boundary
parameter appearing in (2.3), $\varphi(k,x;j)$ is the regular solution
in (2.5), $F(k;j)$ is the Jost function in (2.6), and
$g_j$ is the Gel'fand-Levitan bound-state norming constant in (2.9).

We recall that the $\gamma_s$-values are not necessarily
in an increasing or decreasing order, and the ordering only refers to the order in which
the bound states are added.
We suppose that the bound states
are added in succession by starting with the potential $V(x;0)$ containing
no bound states and by first adding the bound state at $k=i\gamma_1$
with the Gel'fand-Levitan norming constant $g_1,$ then by
adding the bound state at $k=i\gamma_2$ with the norming constant
$g_2,$ and so on. In the presence of
$N$ bound states, when the bound states are removed in succession, we start with
the potential $V(x;N)$ and first remove the bound state at $k=i\gamma_N$
with the norming constant $g_N,$ then remove the bound state at $k=i\gamma_{N-1}$
with the norming constant $g_{N-1},$ and so on.

The following theorem summarizes the Darboux transformation when a bound state
at $k=i\gamma_{j+1}$ with the Gel'fand-Levitan norming constant $g_{j+1}$ is added
to the half-line Schr\"odinger operator with the potential $V(\cdot;j)$
and the boundary parameter $\theta_j.$ In the Dirichlet case, i.e. when $\theta_j=\pi,$
we refer the reader
to [8] for the Darboux transformation formulas provided in the theorem.
In the non-Dirichlet case, i.e. when $\theta_j\in(0,\pi),$
we refer the reader
to (2.3.23) of [15] for the Darboux transformation formulas when a bound state is added.
The formulas in the non-Dirichlet case look similar to those in the
Dirichlet case except that the boundary parameter $\theta_j$ has to be allowed
to change so that the two conditions given in the first line of (2.5)
are satisfied. We invite the interested
reader to directly verify
the results by showing that (2.1) and (2.5) are satisfied
after the bound state is added.

\noindent {\bf Theorem 3.1} {\it Let $V(\cdot;j)$ be the potential of the
Schr\"odinger operator specified in (2.1)-(2.3) with the boundary
parameter $\theta_j$ and the bound states at $k=i\gamma_s$
for $s=1,\dots,j,$ where we assume that there are no bound
states in case $j=0.$
Assume that
one bound state at $k=i\gamma_{j+1}$ is added to the spectrum with
the Gel'fand-Levitan norming constant $g_{j+1},$ but otherwise the relevant
spectral data set is unchanged. The
resulting boundary parameter $\theta_{j+1},$ potential $V(x;j+1),$
regular solution $\varphi(k,x;j+1),$ and Jost
function $F(k;j+1)$ are related to the original quantities $\theta_j,$ $V(x;j),$
$\varphi(k,x;j),$ and $F(k;j)$ as}
$$\cases \cot\theta_{j+1}=\cot\theta_j+g^2_{j+1},&\qquad \theta_j\in(0,\pi),\\
\stretch
\theta_{j+1}=\theta_j,&\qquad \theta_j=\pi,\endcases \tag 3.1$$
$$V(x;j+1)=V(x;j)-\ds\frac{d}{dx}\left[
\ds\frac{2g^2_{j+1}\, \varphi(i\gamma_{j+1},x;j)^2}
{1+g^2_{j+1}\ds\int_0^x dy\,\varphi(i\gamma_{j+1},y;j)^2}
\right],\tag 3.2$$
$$F(k;j+1)=\ds\frac{k-i\gamma_{j+1}}{k+i\gamma_{j+1}}\,F(k;j),
\tag 3.3$$
$$\varphi(k,x;j+1)=\varphi(k,x;j)-
\ds\frac{g^2_{j+1}\, \varphi(i\gamma_{j+1},x;j)\ds\int_0^x dy\,\varphi(k,y;j)\,
\varphi(i\gamma_{j+1},y;j)}{1+g^2_{j+1}\ds\int_0^x dy\,\varphi(i\gamma_{j+1},y;j)^2}.
\tag 3.4$$

The following theorem summarizes the Darboux transformation when the bound state
at $k=i\gamma_j$ with the Gel'fand-Levitan norming constant $g_j$ is removed
 from the half-line Schr\"odinger operator with the potential $V(\cdot;j)$
and the boundary parameter $\theta_j.$
%For the formulas in the
%Dirichlet case, the reader is referred to [8].
The formulas in the
non-Dirichlet case resemble the corresponding formulas in the Dirichlet case
except that the boundary parameter changes in a way compatible with
the first line of (3.1). We omit the proof of the theorem
and invite the interested reader to directly verify the formulas
by showing that (2.1) and (2.5) are satisfied after the bound state is removed.

\noindent {\bf Theorem 3.2} {\it Let $V(\cdot;j)$ be the potential of the
Schr\"odinger operator specified in (2.1)-(2.3) with the boundary
parameter $\theta_j$ and the bound states at
$k=i\gamma_s$ for $s=1,\dots,j.$
Assume that
the bound state at $k=i\gamma_j$ is removed from the spectrum with
the Gel'fand-Levitan norming constant $g_j,$ but otherwise the relevant
spectral data set is unchanged.
The resulting boundary parameter $\theta_{j-1},$
potential $V(x;j-1),$
regular solution $\varphi(k,x;j-1),$ Jost
function $F(k;j-1)$
are related to $\theta_j,$ $V(x;j),$ $\varphi(k,x;j),$ and
$F(k;j)$ as}
$$\cases \cot\theta_{j-1}=\cot\theta_j-g^2_j,&\qquad \theta_j\in(0,\pi),\\
\stretch
\theta_{j-1}=\theta_j,&\qquad \theta_j=\pi,\endcases$$
$$V(x;j-1)=V(x;j)+\ds\frac{d}{dx}\left[
\ds\frac{2g^2_j\, \varphi(i\gamma_j,x;j)^2}{1-g^2_j\ds\int_0^x dy\,\varphi(i\gamma_j,y;j)^2}
\right],\tag 3.5$$
$$F(k;j-1)=\ds\frac{k+i\gamma_j}{k-i\gamma_j}
\,F(k;j),\tag 3.6$$
$$\varphi(k,x;j-1)=\varphi(k,x;j)+
\ds\frac{g^2_j\,\varphi(i\gamma_j,x;j)\ds\int_0^x dy\,\varphi(k,y;j)\,
\varphi(i\gamma_j,y;j)}{1-g^2_j\ds\int_0^x dy\,\varphi(i\gamma_j,y;j)^2},\tag 3.7$$
{\it where $\theta_0,$ $V(x;0),$ $F(k;0),$ and $\varphi(k,x;0)$
correspond to the relevant quantities with no bound
states.}

Let us remark that (2.11), (3.3), and (3.6) imply that the boundary
conditions cannot switch from
a Dirichlet condition to a non-Dirichlet condition or vice versa when bound states are added
or removed via a Darboux transformation. This is because (3.3) and
(3.6) show that the leading term in (2.11) for the large-$k$ asymptotics
of the Jost function $F_\theta(k)$ cannot change from $1$ to $k$ or vice versa
as $k\to+\infty.$

The next theorem indicates that the compact-support property of the potential
is retained if a bound state is removed.

\noindent {\bf Theorem 3.3} {\it Let $V(\cdot;j)\in\Cal A$ be the potential of the
Schr\"odinger operator specified in (2.1)-(2.3) with the boundary
parameter $\theta_j,$ the constant $b$ in (2.2) related
to the compact support of $V(\cdot;j),$ and the bound states at
$k=i\gamma_s$ for $s=1,\dots,j.$
Assume that
the bound state at $k=i\gamma_j$ is removed from the spectrum with
the Gel'fand-Levitan norming constant $g_j,$ but otherwise the relevant
spectral data set is unchanged.
If the compact support of $V(\cdot;j)$ is confined to the interval $(0,b),$ then
the support of $V(\cdot;j-1)$ is also confined to $(0,b)$ and
we have $V(\cdot;j-1)\in\Cal A.$}

\noindent PROOF: We know that (3.5) holds because $V(\cdot;j)$
has a bound state at $k=i\gamma_j$
with the norming constant $g_j$ given in (2.9). It is enough to show that
the quantity inside the brackets in (3.5) is a constant for $x\ge b$
and hence its $x$-derivative vanishes. Because
$\varphi(i\gamma_j,x;j)$ is a bound state, it decays exponentially as
$x\to+\infty.$ Thus, from (2.7), by using (2.4) and
$F(i\gamma_j;j)=0$ we get
$$\varphi(i\gamma_j,x;j)^2=\mp \ds\frac{1}{4\gamma_j^2}\,F(-i\gamma_j;j)^2\,
e^{-2\gamma_j x},\qquad x\ge b,\tag 3.8$$
where the upper sign refers to the non-Dirichlet case $\theta_j\in(0,\pi)$
and the lower
sign to the Dirichlet case $\theta_j=\pi.$
 For $x\ge b$ we can evaluate the denominator
inside the brackets in (3.5) by using $\int_0^x=\int_0^\infty-\int_x^\infty$ there.
Because of (2.9) we have
$$g_j^2 \ds\int_0^\infty dy\,\varphi(i\gamma_j,y;j)^2=1,\tag 3.9$$
and with the help of (3.8) we get
$$g_j^2 \ds\int_x^\infty dy\,\varphi(i\gamma_j,y;j)^2=\mp \ds\frac{1}{8\gamma_j^3}\,F(-i\gamma_j;j)^2\,
e^{-2\gamma_j x},\qquad x\ge b.\tag 3.10$$
Using (3.8)-(3.10) in the quantity inside the brackets in (3.5), we get
$$
\ds\frac{2g^2_j\, \varphi(i\gamma_j,x;j)^2}{1-g^2_j\ds\int_0^x dy\,\varphi(i\gamma_j,y;j)^2}=4\gamma_j,\qquad x\ge b,\tag 3.11$$
and hence from (3.5) we see that $V(x;j)\equiv V(x;j-1)$ for $x> b$
and thus $V(\cdot;j-1)$ has the same support as $V(\cdot;j).$
The property
$V(\cdot;j-1)\in\Cal A$ then follows from
the fact that the quantity inside the brackets
in the second term on the right-hand side of (3.5)
is real valued and continuous in $x$ when $x\in[0,b].$ \qed

In the notation used in this section, we can express the definition of
the Gel'fand-Levitan norming constant $g_s$ given in (2.9) as
$$g_s:=\ds\frac{1}{\ds\sqrt{\int_0^\infty dx\,\varphi(i\gamma_s,x;N)^2}},
\qquad s=1,\dots,N,\tag 3.12$$
where $\varphi(k,x;N)$ is the regular solution appearing in (2.7).
The following result shows that we can obtain $g_j$ by normalizing
not only $\varphi(i\gamma_j,x;N)$ but any one of
$\varphi(i\gamma_j,x;s)$ for $s=j,j+1,\dots,N.$

\noindent {\bf Theorem 3.4} {\it Let $V(\cdot;N)\in\Cal A$ be the potential of the
Schr\"odinger operator specified in (2.1)-(2.3) with the boundary
parameter $\theta_N,$ the bound states at
$k=i\gamma_s$ for $s=1,\dots,N,$ and
the corresponding Gel'fand-Levitan norming
constants $g_s$ defined as in (2.9). For any $j$ with $1\le j< N,$ we then have}
$$\ds\int_0^\infty dx\, \varphi(i\gamma_j,x;j)^2=\ds\int_0^\infty dx\, \varphi(i\gamma_j,x;j+1)^2
=\cdots=
\ds\int_0^\infty dx\, \varphi(i\gamma_j,x;N)^2.\tag 3.13$$

\noindent PROOF: From (3.4), for any positive integer $s$
with $j+1\le s\le N,$ we obtain
$$\varphi(i\gamma_j,x;s)=\varphi(i\gamma_j,x;j)-
\ds\frac{g_s^2\,\varphi(i\gamma_s,x;j)\,\ds\int_0^x dy\,
\varphi(i\gamma_j,y;j)\,\varphi(i\gamma_s,y;j)}{1+g^2_s\ds\int_0^x dy\,\varphi(i\gamma_s,y;j)^2}.\tag 3.14$$
Squaring both sides of (3.14) and with some simplification,
we observe that
$$\varphi(i\gamma_j,x;s)^2=\varphi(i\gamma_j,x;j)^2-\ds\frac{d}{dx}\left[
\ds\frac{g_s^2\,\left[\ds\int_0^x dy\,\varphi(i\gamma_j,y;j)\,
\varphi(i\gamma_s,y;j)\right]^2}{1+g^2_s\ds\int_0^x dy\,\varphi(i\gamma_s,y;j)^2}\right].\tag 3.15$$
Integrating both sides of (3.15) over $x\in(0,+\infty),$
we see that the equalities in (3.13) all hold provided the quantity
inside the brackets in (3.15) vanishes as $x\to+\infty$
because that quantity already vanishes at $x=0.$
Let us use $\int_0^x=\int_0^b+\int_b^x$ when $x\ge b$ and
estimate the integrals in the numerator and in the denominator in (3.15).
By Theorem~3.3 we know that $V(\cdot;j)\in \Cal A$ because
$V(\cdot;N)\in \Cal A.$ Thus, $V(x;j)\equiv 0$ for $x>b$ and $f(k,x;j)=e^{ikx}$
for $x\ge b$ as a result of (2.4). We also have
$F(i\gamma_j;j)=0$ and thus via (3.3) we have $F(i\gamma_s;j)\ne 0$ for $j+1
\le s\le N.$
Therefore, from (2.7) we obtain
$$\varphi(i\gamma_j,x;j)^2=\mp\ds\frac{1}{4\gamma_j^2}\,
F(-i\gamma_j;j)^2\,e^{-2\gamma_j x},\qquad x\ge b,
\tag 3.16$$
and for $j+1\le s\le N$ we have
$$\varphi(i\gamma_s,x;j)^2=\mp\ds\frac{1}{4\gamma_s^2}\,
\left[F(i\gamma_s;j)\,e^{\gamma_s x}-F(-i\gamma_s;j)\,e^{-\gamma_s x}\right]^2,\qquad x\ge b,
\tag 3.17$$
where the upper sign refers to the non-Dirichlet case $\theta_j\in (0,\pi)$ and
the lower sign refers to the Dirichlet case $\theta_j=\pi.$
With the help of (3.16) and (3.17) we get
$$\ds\int_0^x dy\,\varphi(i\gamma_j,y;j)\,
\varphi(i\gamma_s,y;j)=O\left(e^{(\gamma_s-\gamma_j) x}\right),\qquad x\to +\infty,$$
$$\ds\int_0^x dy\,\varphi(i\gamma_s,y;j)^2
=O\left(e^{2 \gamma_s  x}\right),\qquad x\to +\infty.$$
Thus, the quantity inside the brackets in (3.15) has the behavior
$O(e^{-2 \gamma_j  x})$ as $x\to+\infty.$ Hence, our proof is complete. \qed

Using the result in Theorem~3.4 we can comment on the denominator in
(3.5). As seen from (3.12) and (3.13), the Gel'fand-Levitan norming constant
$g_j$ can be
obtained by normalizing $\varphi(i\gamma_j,x;s)$ for any
integer $s$ with $j\le s\le N,$ i.e. via
$$g_j=\ds\frac{1}{\ds\sqrt{\int_0^\infty dx\,\varphi(i\gamma_j,x;s)^2}},
\qquad
s=j,j+1,\dots,N.\tag 3.18$$
Using (3.5) and the positivity of
$\varphi(i\gamma_j,y;j)^2$, we conclude that
the integral $\int_0^x dy\,\varphi(i\gamma_j,y;j)^2$ is an
increasing function of $x.$ With the help of (3.18) we see that
it increases from the value of zero at $x=0$ to the
value of $1/g_j^ 2$ as $x$ increases from $x=0$ to $x=+\infty.$
Thus, the denominator in (3.5) remains positive for $x\in\bR^+.$

The following theorem is one of the key results needed
for the characterization of eligible and ineligible resonances.
Recall that an eligible resonance corresponds to a zero of the Jost function
defined in (2.6) in such a way that such a zero occurs on the negative
imaginary axis and
can be converted into a
bound state through a Darboux transformation without changing the
compact support of the potential. If a zero of
the Jost function occurring on
the negative imaginary axis cannot be converted into a bound state under a Darboux transformation without changing the
compact support of the potential, then we refer to such an imaginary resonance
as an ineligible resonance.

\noindent {\bf Theorem 3.5} {\it Let $V(\cdot;j)\in\Cal A$ be the potential of the
Schr\"odinger operator specified in (2.1)-(2.3) with the boundary
parameter $\theta_j$ and the bound states at
$k=i\gamma_s$ for $s=1,\dots,j.$
Assume that
a bound state at $k=i\gamma_{j+1}$
is added to the spectrum with
the Gel'fand-Levitan norming constant $g_{j+1},$ but otherwise the relevant
spectral data set is unchanged. Let $b$ be the constant appearing
in (2.2) related to the compact support of $V(\cdot;j).$
The support of $V(\cdot;j+1)$ is also confined to $(0,b)$ if and only if}
$$F(-i\gamma_{j+1};j)=0,\quad
g_{j+1}^2=\ds\frac{2\gamma_{j+1}}{\varphi(i\gamma_{j+1},b;j)^2-2\gamma_{j+1}\ds\int_0^b dy\,
\varphi(i\gamma_{j+1},y;j)^2}.\tag 3.19$$
{\it Note that the second condition in (3.19) implies that we must have}
$$\ds\frac{2\gamma_{j+1}}{\varphi(i\gamma_{j+1},b;j)^2-2\gamma_{j+1}\ds\int_0^b dy\,
\varphi(i\gamma_{j+1},y;j)^2}>0.\tag 3.20$$
{\it When (3.19) is satisfied, the resulting potential $V(\cdot;j+1)$ belongs to class $\Cal A.$}

\noindent PROOF: In order to prove our theorem, from (3.2) we see that it is enough
to prove that (3.19) is equivalent to $F(i\gamma_{j+1};j+1)=0$ and that
$$\ds\frac{2g_{j+1}^2 \, \varphi
(i\gamma_{j+1},x;j)^2}{1+g_{j+1}^2\ds\int_0^x dy\,
\varphi(i\gamma_{j+1},y;j)^2}=c_1,
\qquad x\ge b,\tag 3.21$$
for some constant $c_1.$
In fact, from (3.11) we know that the value of $c_1$ must be $4\gamma_{j+1}.$
We first show that (3.19) holds with $c_1=4\gamma_{j+1}$ there. For this we proceed as follows.
Because $k=i\gamma_{j+1}$ corresponds to a bound state, we have $F(i\gamma_{j+1};j+1)=0.$ By Theorem~2.1 we know that
$F(k;j)$ is entire in $k,$ and hence
 from (3.3) we see that we must have $F(-i\gamma_{j+1};j)=0.$
Since $V(x;j)\equiv 0$ for $x>b,$ by (2.4) the corresponding Jost solution is given by
$f(k,x;j)=e^{ikx}$ for $x\ge b.$ Using $F(-i\gamma_{j+1};j)=0$ in (2.7) we get
$$\varphi(i\gamma_{j+1},x;j)^2=\mp \ds\frac{1}{4\gamma_{j+1}^2}\,F(i\gamma_{j+1};j)^2\,
e^{-2\gamma_{j+1} x},\qquad x\ge b,$$
where the upper sign refers to the non-Dirichlet case
$\theta_j\in(0,\pi)$ and the lower
sign to the Dirichlet case $\theta_j=\pi.$
Thus, (3.21) is satisfied provided we have
$$\aligned 4 & \gamma_{j+1}=\\ &
\ds\frac{\mp\ds\frac{g_{j+1}^2}{2\gamma_{j+1}^2}\, F(i\gamma_{j+1};j)^2
\,e^{-2\gamma_{j+1} x}}
{1+g_{j+1}^2 \ds\int_0^b dy\, \varphi(i\gamma_{j+1},y;j)^2
-\ds\frac{g_{j+1}^2}{2\gamma_{j+1}}\,\varphi(i\gamma_{j+1},b;j)^2
\mp \frac {g_{j+1}^2}
{8\gamma_{j+1}^3}\, F(i\gamma_{j+1};j)^2\,e^{-2\gamma_{j+1} x} }.
\endaligned\tag 3.22$$
After cross multiplying and simplifying, we see that (3.22) is equivalent to
$$1+g_{j+1}^2 \ds\int_0^b dy\, \varphi(i\gamma_{j+1},y;j)^2
-\ds\frac{g_{j+1}^2}{2\gamma_{j+1}}\,\varphi(i\gamma_{j+1},b;j)^2
=0,$$
which is satisfied because of the second equality in (3.19).
Let us now prove the converse, namely, prove that
 $V(x;j+1)\equiv 0$ for $x>b$
implies (3.19). From (3.2) and the fact that $V(\cdot;j)\in\Cal A$ we know that
 $V(x;j+1)\equiv 0$ for $x>b$ if and only if
(3.21) holds with $c_1=4\gamma_{j+1}$ there, i.e.
$$\ds\frac{2g_{j+1}^2 \, \varphi
(i\gamma_{j+1},x;j)^2}{1+g_{j+1}^2\ds\int_0^x dy\,
\varphi(i\gamma_{j+1},y;j)^2}=4\gamma_{j+1}
,\qquad x\ge b.\tag 3.23$$
Evaluating (3.23) at $x=b$ we get the second equality in (3.19).
Let us cross multiply in (3.23) and then take the $x$-derivative
of both sides of the resulting equation. We get
$$4 g_{j+1}^2 \, \varphi'
(i\gamma_{j+1},x;j)\, \varphi
(i\gamma_{j+1},x;j)=4 g_{j+1}^2 \,\gamma_{j+1}\, \varphi
(i\gamma_{j+1},x;j)^2,\qquad x\ge b,$$
or equivalently
$$\varphi'
(i\gamma_{j+1},x;j)=\gamma_{j+1}\,\varphi
(i\gamma_{j+1},x;j),\qquad x\ge b.\tag 3.24$$
 From (3.24) we see that
 $$\varphi'
(i\gamma_{j+1},x;j)=c_2\, e^{\gamma_{j+1} x},\qquad x\ge b,\tag 3.25$$
for some constant $c_2.$
On the other hand, with the help of (2.4) and (2.7) we get
for $x\ge b$
$$\varphi(i\gamma_{j+1},x;j)=\cases
\ds\frac{1}{2i\gamma_{j+1}}\left[
F(i\gamma_{j+1};j)\,e^{\gamma_{j+1}x}-F(-i\gamma_{j+1};j)
\,e^{-\gamma_{j+1}x}\right],&\qquad \theta\in
(0,\pi),\\
\stretch
\ds\frac{1}{2\gamma_{j+1}}\left[
F(i\gamma_{j+1};j)\,e^{\gamma_{j+1}x}-F(-i\gamma_{j+1};j)
\,e^{-\gamma_{j+1}x}\right],&\qquad \theta=
\pi.\endcases\tag 3.26$$
Comparing (3.25) and (3.26) we see that we must have $F(-i\gamma_{j+1};j)=0.$ When (3.22)
is satisfied, the potential $V(\cdot;j+1)$ belongs to $\Cal A$ because
the quantity inside the brackets
in the second term on the right-hand side of (3.2)
is real valued and continuous in $x$ when $x\in[0,b].$
 \qed

The result in the preceding theorem is fascinating in the sense that if we add a bound state to the compactly-supported potential $V(\cdot;j)$ in class $\Cal A$
at $k=i\gamma_{j+1}$ with some arbitrary Gel'fand-Levitan norming constant $g_{j+1},$
in general the resulting
potential $V(\cdot;j+1)$ cannot be compactly supported. Theorem~3.5 states that the potential $V(\cdot;j+1)$
is compactly supported if and only if
$k=-i\gamma_{j+1}$ happens to be a zero of $F(k;j)$
and the norming constant $g_{j+1}$ happens to be equal to the square root of the
quantity on the right-hand side of the second
equality in (3.19).
Thus, if the left-hand side in (3.20) does not yield
a positive
number, then it is impossible for $V(\cdot;j+1)$ to have the support in $(0,b)$
because there cannot be a corresponding positive norming constant $g_{j+1}$ guaranteeing
the compact support for the potential.
Let us clarify that, if the left-hand side in (3.20) is not positive, one can find a potential
with support in $(0,b),$ but such a potential must have a singularity and it cannot belong
to class $\Cal A.$

The result of Theorem~3.5 is
analogous to the result [9] from the full-line Schr\"odinger
equation when a bound state is added to a compactly-supported potential: Start with a
compactly-supported potential $V$ associated with
the transmission coefficient $T$
and add a bound state to it at $k=i\kappa$
to obtain the potential $\tilde V$ with the transmission coefficient
$\tilde T$ given by
$$\tilde T(k)=\ds\frac{k+i\kappa}{k-i\kappa}\,T(k).$$
Then, $\tilde V$ is also compactly
supported if and only if the transmission coefficient $T(k)$
has a pole at $k=-i\kappa.$ The analysis in the full-line case is less complicated
due to the fact that
in the full-line case there is no boundary condition at $x=0$ such as (2.3).

In Theorem~3.5, in terms of $F(k;j)$ and $\varphi(k,x;j),$ we have expressed the necessary and sufficient conditions
for the potential $V(\cdot;j+1)$ to have the same compact support as
$V(\cdot;j).$ In the next theorem the two conditions stated in
(3.19) are expressed in terms of $F(k;j+1)$ and $\varphi(k,x;j+1).$

\noindent {\bf Theorem 3.6} {\it Let $V(\cdot;j)\in\Cal A$ be the potential of the
Schr\"odinger operator specified in (2.1)-(2.3) with the boundary
parameter $\theta_j$ and the bound states at
$k=i\gamma_s$ for $s=1,\dots,j.$
Assume that
a bound state at $k=i\gamma_{j+1}$
is added to the spectrum with
the Gel'fand-Levitan norming constant $g_{j+1},$ but otherwise the relevant
spectral data set is unchanged. Let $b$ be the constant appearing
in (2.2) related to the compact support of $V(\cdot;j).$
The support of $V(\cdot;j+1)$ is also confined to $(0,b)$ if and only if}
$$F(i\gamma_{j+1};j+1)=0,\quad
g_{j+1}^2=\ds\frac{2\gamma_{j+1}}{\varphi(i\gamma_{j+1},b;j+1)^2+2\gamma_{j+1}\ds\int_0^b dt\,
\varphi(i\gamma_{j+1},t;j+1)^2}.\tag 3.27$$

\noindent PROOF: The equivalence of $F(i\gamma_{j+1};j+1)=0$ and $F(-i\gamma_{j+1};j)=0$ is already shown in the
proof of Theorem~3.5.
Let us now prove that the second equality in (3.19) is equivalent to
the second equality in (3.27). From (3.5) we see that
$$V(x;j+1)=V(x;j)-\ds\frac{d}{dx}\left[
\ds\frac{2g_{j+1}^2 \, \varphi
(i\gamma_{j+1},x;j+1)^2}{1-g_{j+1}^2\ds\int_0^x dy\,
\varphi(i\gamma_{j+1},y;j+1)^2}
\right].\tag 3.28$$
A comparison with (3.2) shows that the right-hand sides of
(3.2) and of (3.28) are equal to each other for $x>b,$ and we have
$$
\ds\frac{2g_{j+1}^2 \, \varphi
(i\gamma_{j+1},x;j+1)^2}{1-g_{j+1}^2\ds\int_0^x dy\,
\varphi(i\gamma_{j+1},y;j+1)^2}=
\ds\frac{2g_{j+1}^2 \, \varphi
(i\gamma_{j+1},x;j)^2}{1+g_{j+1}^2\ds\int_0^x dy\,
\varphi(i\gamma_{j+1},y;j)^2}+c_3,\qquad x\ge b,\tag 3.29$$
for some constant $c_3.$ Using (3.23) on the right-hand side of (3.29),
we get
$$
\ds\frac{2g_{j+1}^2 \, \varphi
(i\gamma_{j+1},x;j+1)^2}{1-g_{j+1}^2\ds\int_0^x dy\,
\varphi(i\gamma_{j+1},y;j+1)^2}=
4\gamma_{j+1}+c_3,\qquad x\ge b.\tag 3.30$$
Cross multiplying in (3.30) and then taking the $x$-derivative
of the resulting equation, for $x\ge b$ we obtain
$$4\,g_{j+1}^2 \, \varphi'
(i\gamma_{j+1},x;j+1)\, \varphi
(i\gamma_{j+1},x;j+1)=-(4\gamma_{j+1}+c_3)g_{j+1}^2\,
\varphi
(i\gamma_{j+1},x;j+1)^2,$$
which simplifies to
$$\varphi'
(i\gamma_{j+1},x;j+1)=-\left(\gamma_{j+1}+\ds\frac{c_3}{4}\right)
\varphi'
(i\gamma_{j+1},x;j+1),\qquad x\ge b.\tag 3.31$$
On the other hand, since $V(x;j+1)\equiv 0$ for $x>b,$
we have the analog of (3.16) given by
$$\varphi(i\gamma_{j+1},x;j+1)^2=\mp\ds\frac{1}{4\gamma_{j+1}^2}\,
F(-i\gamma_{j+1};j)^2\,e^{-2\gamma_{j+1} x},\qquad x\ge b,
\tag 3.32$$
where the upper sign refers to the non-Dirichlet case $\theta_{j+1}\in (0,\pi)$ and
the lower sign to the Dirichlet case $\theta_{j+1}=\pi.$
Comparing (3.31) and (3.32) we get $c_3=0,$ and hence
(3.30) yields
$$\ds\frac{2g_{j+1}^2 \, \varphi
(i\gamma_{j+1},b;j+1)^2}{1-g_{j+1}^2\ds\int_0^b dy\,
\varphi(i\gamma_{j+1},y;j+1)^2}=
4\gamma_{j+1}.\tag 3.33$$
By isolating $g_{j+1}^2$ to one side of the equation in (3.33),
we observe from (3.29) and (3.33) that
the second equality in (3.19) is equivalent to the second equality in (3.27). \qed

We can ask whether we can predict if (3.20) is satisfied without actually
evaluating the left-hand side in (3.20).
For this purpose, we will exploit the signs of $\varphi(i\gamma_{j+1},x;j)$
and $\varphi(i\gamma_{j+1},x;j+1)$ as $x\to+\infty.$ It is convenient
to define
$$H(\beta;j):=\cases -i\,F(i\beta;j),&\qquad
\theta_j\in(0,\pi),\\
\stretch
F(i\beta;j),&\qquad
\theta_j=\pi,\endcases\tag 3.34$$
where $F(k;j)$ is the Jost function corresponding to the potential
$V(\cdot;j)$ and the boundary parameter $\theta_j.$
The advantage of using $H(\beta;j)$ rather than $F(i\beta;j)$
is that the former is real valued and hence its sign can be examined graphically.
Note that
$$H'(\beta;j):=\ds\frac{dH(\beta;j)}{d\beta}
=\cases
\left.\ds\frac{dF(k;j)}{dk}\ds\right|_{k=i\beta},
&\qquad \theta_j\in(0,\pi),\\
i\left.\ds\frac{dF(k;j)}{dk}\ds\right|_{k=i\beta},
&\qquad \theta_j=\pi.\endcases\tag 3.35$$
Note also that,
as seen from (2.25), as $\beta\to+\infty$ we have
$$H(\beta;j)=\cases
\beta+O(1),&\qquad \theta_j\in(0,\pi),\\
\stretch
1+O\left(\ds\frac{1}{\beta}\right),&\qquad
\theta_j=\pi,\endcases\tag 3.36$$
and hence $H(\beta;j)$ is positive for large
positive $\beta$-values.

The result in the following theorem can be used as a test to determine
whether the inequality in (3.20) is satisfied or not.

\noindent {\bf Theorem 3.8} {\it Let $V(\cdot;j)\in\Cal A$ be the potential of the
Schr\"odinger operator specified in (2.1)-(2.3) with the boundary
parameter $\theta_j$ and the bound states at
$k=i\gamma_s$ for $s=1,\dots,j.$ Let $F(k;j)$ be the corresponding Jost
function defined in (2.6), $H(\beta;j)$ be the quantity
defined in (3.34), and $b$ be the constant appearing in (2.2).
Assume that
a bound state at $k=i\gamma_{j+1}$
is added to the spectrum,
but otherwise the relevant
spectral data set is unchanged.
The support of $V(\cdot;j+1)$ is also confined to the interval
$(0,b)$ if and only if}
$$F(-i\gamma_{j+1};j)=0,\quad \ds\frac{i\,\dot F(-i\gamma_{j+1};j)}{F(i\gamma_{j+1};j)}>0,\tag 3.37$$
{\it or equivalently, if and only if}
$$H(-\gamma_{j+1};j)=0,\quad \ds\frac{H'(-\gamma_{j+1};j)}
{H(\gamma_{j+1};j)}>0.\tag 3.38$$

\noindent PROOF: The equivalence of (3.37) and (3.38) is obtained directly by using (3.34) and (3.35).
Thus, we only need to show
that (3.37) is equivalent to the first condition given in (3.19) and
the condition in (3.20). In other words,
we need to prove that (3.37) is equivalent to
$$F(-i\gamma_{j+1};j)=0,\tag 3.39$$
and to the positivity of the
right-hand side in the equality involving
$g_{j+1}^2$ in (3.19).
Note that (3.39) appears also in (3.19) and hence we only need to show the equivalence of
the inequality in (3.37) and the positivity of the relevant
quantity. Using (3.39) in
(3.26) we see that, for $x\ge b,$ we have
$$\varphi(i\gamma_{j+1},x;j)=\cases
\ds\frac{1}{2i\gamma_{j+1}}\,
F(i\gamma_{j+1};j)\,e^{\gamma_{j+1}x},&\qquad \theta_j\in (0,\pi)
,\\
\stretch
\ds\frac{1}{2\gamma_{j+1}}\,
F(i\gamma_{j+1};j)\,e^{\gamma_{j+1}x},&\qquad \theta_j=
\pi.\endcases\tag 3.40$$
Using (3.39) we can write (3.3) as
$$F(k,j+1)=(k-i\gamma_{j+1})\,\ds\frac{F(k;j)-F(-i\gamma_{j+1};j)}{k+i\gamma_{j+1}}.\tag 3.41$$
Letting $k\to -i\gamma_{j+1},$ from (3.41), as a result of the
analyticity of $F(k;j)$ in $\bC$ we obtain
$$F(-i\gamma_{j+1},j+1)=-2i\gamma_{j+1}\,\dot F(-i\gamma_{j+1};j),\tag 3.42$$
where we recall that an overdot indicates the $k$-derivative.
With the help of (2.7)
let us now evaluate $\varphi(i\gamma_{j+1},x;j+1).$ Using (2.4) in (2.7), for $x\ge b$ we obtain
$$\varphi(i\gamma_{j+1},x;j+1)=\cases
\ds\frac{1}{2i\gamma_{j+1}}\left[
F(i\gamma_{j+1};j+1)\,e^{\gamma_{j+1}x}-F(-i\gamma_{j+1};j+1)
\,e^{-\gamma_{j+1}x}\right],\\
\stretch
\ds\frac{1}{2\gamma_{j+1}}\left[
F(i\gamma_{j+1};j)\,e^{\gamma_{j+1}x}-F(-i\gamma_{j+1};j+1)
\,e^{-\gamma_{j+1}x}\right],\endcases\tag 3.43$$
where the first line holds if $\theta_{j+1}\in (0,\pi)$ and the second
line holds if $\theta_{j+1}=\pi.$
 From Theorem~3.6 we know that $F(i\gamma_{j+1},j+1)=0$ and hence (3.43), for $x\ge b,$ is equivalent to
$$\varphi(i\gamma_{j+1},x;j+1)=\cases
-\ds\frac{1}{2i\gamma_{j+1}}\,
F(-i\gamma_{j+1};j+1)
\,e^{-\gamma_{j+1}x},&\qquad \theta_j\in (0,\pi),\\
\stretch
-\ds\frac{1}{2\gamma_{j+1}}\,
F(-i\gamma_{j+1};j+1)
\,e^{-\gamma_{j+1}x},&\qquad \theta_j=
\pi.\endcases\tag 3.44$$
Using (3.42) in (3.44) we see that, for $x\ge b,$ we have
$$\varphi(i\gamma_{j+1},x;j+1)=\cases
\dot F(-i\gamma_{j+1};j)
\,e^{-\gamma_{j+1}x},&\qquad \theta_j\in (0,\pi),\\
\stretch
i\,\dot
F(-i\gamma_{j+1};j)
\,e^{-\gamma_{j+1}x},&\qquad \theta_j=
\pi.\endcases\tag 3.45$$
With the help of (3.1), we see that $\theta_{j+1}\in (0,\pi)$ if and only if $\theta_j\in (0,\pi).$
Hence, from (3.40) and (3.45) we obtain
$$\ds\frac{\varphi(i\gamma_{j+1},x;j+1)}{\varphi(i\gamma_{j+1},x;j)}=
2i\gamma_{j+1}\,\ds\frac{\dot F(-i\gamma_{j+1};j)}{F(i\gamma_{j+1};j)}\,e^{-\gamma_{j+1}x},
\qquad \theta_j\in(0,\pi],\quad x\ge b.\tag 3.46$$
 From (3.46) we see that the inequality in (3.37) is satisfied if and only if the quantity on the
 left-hand side of (3.46) is positive for any $x\ge b.$ Let us now evaluate that
 quantity. From (3.7), using $j+1$ instead of $j$ there and letting $k=i\gamma_{j+1}$
 there, we obtain
 $$\varphi(i\gamma_{j+1},x;j)=
\ds\frac{\varphi
(i\gamma_{j+1},x;j+1)}{1-g_{j+1}^2\ds\int_0^x dy\,
\varphi(i\gamma_{j+1},y;j+1)^2},\qquad x\ge 0,$$
or equivalently
$$\ds\frac{\varphi
(i\gamma_{j+1},x;j+1)}{\varphi(i\gamma_{j+1},x;j)}=
1-g_{j+1}^2\ds\int_0^x dy\,
\varphi(i\gamma_{j+1},y;j+1)^2,\qquad x\ge 0.\tag 3.47$$
 From (3.18) it follows that
$$\ds\frac{1}{g_{j+1}^2}=\int_0^\infty dy\,
\varphi(i\gamma_{j+1},y;j+1)^2,\tag 3.48$$
and hence using $\int_0^x=\int_0^\infty-\int_x^\infty$ in (3.47), with the help of
(3.48) we get
$$\ds\frac{\varphi
(i\gamma_{j+1},x;j+1)}{\varphi(i\gamma_{j+1},x;j)}=
g_{j+1}^2\ds\int_x^\infty dy\,
\varphi(i\gamma_{j+1},y;j+1)^2,\qquad x\ge 0.\tag 3.49$$
Comparing (3.49) with (3.46) we see that the inequality in (3.37)
is satisfied if and only if $g_{j+1}^2$ appearing in (3.49) is positive.
 From (3.19) and (3.20) we already know that (3.39) and the positivity
 of $g_{j+1}^2$ are equivalent for having $V(\cdot;j+1)$ to have support in $(0,b).$
 Thus, we have proved that (3.19) is equivalent to
 $$F(-i\gamma_{j+1};j)=0; \quad \ds\frac{\varphi
(i\gamma_{j+1},x;j+1)}{\varphi(i\gamma_{j+1},x;j)}>0,\qquad x\ge b.\tag 3.50$$
With the help of (3.46), we see that (3.50) is equivalent to (3.37). Thus, the proof
is complete. \qed

One consequence of Theorem~3.8 is that the scattering matrix corresponding to
a half-line Schr\"odinger operator has a meromorphic extension with simple
poles at the bound states.

\noindent {\bf Proposition 3.9} {\it Let $V(\cdot;j)\in\Cal A$ be the potential of the
Schr\"odinger operator specified in (2.1)-(2.3) with the boundary
parameter $\theta_j$ and the bound states at
$k=i\gamma_s$ for $s=1,\dots,j.$ Let $F(k;j)$ and $S(k;j)$
be the corresponding Jost
function and the scattering matrix defined in (2.6)
and (2.16),respectively.
Assume that a bound state at $k=i\gamma_{j+1}$ is added
to the spectrum without changing the support of the potential
and without changing the remaining part of the spectral data set.
Under the corresponding Darboux transformation, the scattering matrix
is transformed as}
$$S(k;j+1)=\left(\ds\frac{k+i\gamma_{j+1}}{k-i\gamma_{j+1}}\right)^2 S(k;j)
.\tag 3.51$$
{\it The scattering matrix $S(k;j)$ has a meromorphic extension
 from $k\in\bR$ to the entire complex plane. The only poles
 of $S(k;j)$ in $\bCp$ occur at the bound states at
$k=i\gamma_s$ for $s=1,\dots,j$ and such poles are all simple.
Furthermore, $S(k;j)$ has simple zeros at $k=-i\gamma_s$ for $s=1,\dots,j.$}

\noindent PROOF: The meromorphic
extension of $S(k;j)$ from $k\in\bR$ to $k\in\bC$
has already been established in Theorem~2.1(e).
We get (3.51) by using (3.3) in (2.16).
Using induction, from (3.51) it is seen that it is enough to prove that
$S(k;0)$ has no poles in $\bCp$ and that
$S(k;j+1)$ has a simple pole at $k=i\gamma_{j+1}$ and has a simple zero
at $k=-i\gamma_{j+1}.$ Note that
$S(k;0)$ has no poles in $\bCp,$ which follows from (2.16)
and the fact that $F(k;0)$ has no zeros in $\bCp.$
At first sight,
(3.51) gives the wrong impression that $S(k;j+1)$ has a double pole at $k=i\gamma_{j+1}$ and a double zero at $k=-i\gamma_{j+1}.$ However, the
pole at $k=i\gamma_{j+1}$ is a simple one and the zero at $k=-i\gamma_{j+1}$ is a simple one, as the following argument
shows. Using (2.16), let us write (3.51) as
$$S(k;j+1)=\mp\left(\ds\frac{k+i\gamma_{j+1}}{k-i\gamma_{j+1}}\right)\left(
\ds\frac{ F(-k;j)}{k-i\gamma_{j+1}}\right)
\left(
\ds\frac{k+i\gamma_{j+1}}{ F(k;j)}\right)
,\tag 3.52$$
where the upper sign refers to the non-Dirichlet boundary condition
$\theta_j\in(0,\pi)$
and the lower sign to the Dirichlet boundary condition $\theta_j=\pi.$
 From (3.37) we know that $F(-k;j)$ has a simple zero at $k=i\gamma_{j+1}.$
Thus, the second factor on the right-hand side of (3.52)
has a removable singularity at $k=i\gamma_{j+1}$ and no zero
at $k=i\gamma_{j+1}.$ Similarly, the third factor
on the right-hand side of (3.52)
has a removable singularity at $k=-i\gamma_{j+1}$ and no zero
at $k=-i\gamma_{j+1}.$
 We also know that
$F(k;j)$ in the third factor cannot vanish at $k=i\gamma_{j+1}$
 because we already have $F(-i\gamma_{j+1};j)=0$
 as a result of the fact that $F(-k;j)$ and $F(k;j)$ cannot
 vanish at the same $k$-value.
 Thus, the product of the second and third factors on the
 right-hand side of (3.52) does not have a pole at $k=i\gamma_{j+1}$
 and that product does not have a zero at $k=-i\gamma_{j+1}.$
 Hence, the simple pole at $k=i\gamma_{j+1}$ in the first
 factor on the
 right-hand side of (3.52) is the only pole of $S(k;j+1)$ at $k=i\gamma_{j+1}$ and
 that the simple zero at $k=-i\gamma_{j+1}$ in the first
 factor is the only zero of $S(k;j+1)$ at $k=-i\gamma_{j+1}.$
%
% Thus, the second factor on the right-hand side of (3.52)
% has a removable singularity at $k=i\gamma_{j+1}$ and no zero
% there
% and the third
% factor there
%
%
% , neither of
%the second and third factors on the
%right-hand side of (3.52) has a zero at $k=-i\gamma_{j+1},$
%and neither has a pole at $k=-i\gamma_{j+1}.$
%Thus, $S(k;j+1)$ has only a simple pole at $k=i\gamma_{j+1}$ and
%a simple zero at $k=-i\gamma_{j+1}.$ By Theorem~2.1(a) we already know
%that $F(k;j)$ is entire and its only zeros in $\bCpb\setminus\{0\}$
%occur at $k=i\gamma_s$ for $s=1,\dots,j.$ Thus, from (3.52) it
%follows that the only poles
% of $S(k;j)$ in $\bCp$ occur at
%$k=i\gamma_s$ for $s=1,\dots,j$ and such poles are all simple.
%Since the second factor on the right-hand side
% of (3.52) has simple zeros at $k=-i\gamma_s$ for $s=1,\dots,j$
% and $F(k;j)$ has no poles in $\bC,$ it follows from (3.52) that the zeros
% of $S(k;j)$ at $k=-i\gamma_s$ for $s=1,\dots,j$ are all simple.
 \qed

It is useful to state the result of Theorem~3.8 in terms of the quantities associated
with no bound states. Thus, we present the following result.

\noindent {\bf Theorem 3.10} {\it Let $V(\cdot;j)\in\Cal A$ be the potential of the
Schr\"odinger operator specified in (2.1)-(2.3) with the boundary
parameter $\theta_j$ and the bound states at
$k=i\gamma_s$ for $s=1,\dots,j.$ Let $H(\beta;j)$ be the quantity
defined in (3.34) and $b$ be the constant appearing in (2.2).
Assume that
a bound state at $k=i\gamma_{j+1}$
is added to the spectrum,
but otherwise the
spectral data set is unchanged.
The support of $V(\cdot;j+1)$ is also confined to $(0,b)$ if and only if}
$$H(-\gamma_{j+1};0)=0,\quad H'(-\gamma_{j+1};0)>0,\tag 3.53$$
where we recall that $H(\beta;0)$ refers to the quantity in
(3.34) when there are no bound states.

\noindent PROOF: From (3.3) and (3.34) we obtain
$$H(\beta;j)=H(\beta;0) \ds\prod_{s=1}^j\left( \ds\frac{\beta-\gamma_s}{\beta+\gamma_s}
\right).
\tag 3.54$$
Thus, through differentiation with respect to $\beta,$ (3.54) yields
$$H'(\beta;j)=H(\beta;j)\ds\sum_{s=1}^j\left(\ds\frac{2\gamma_s}{\beta+\gamma_s}\right)+
H'(\beta;0)\ds\prod_{s=1}^j\left( \ds\frac{\beta-\gamma_s}{\beta+\gamma_s}
\right).\tag 3.55$$
 From (3.54) and (3.55) we obtain
$$H(\gamma_{j+1};j)=H(\gamma_{j+1};0) \ds\prod_{s=1}^j\left( \ds\frac{\gamma_{j+1}-\gamma_s}{\gamma_{j+1}+\gamma_s}
\right),
\tag 3.56$$
$$H'(-\gamma_{j+1};j)=
H'(-\gamma_{j+1};0)\ds\prod_{s=1}^j\left( \ds\frac{\gamma_{j+1}+\gamma_s}{\gamma_{j+1}-\gamma_s}
\right),\tag 3.57$$
where we have used $H(-\gamma_{j+1};j)=0$ to get (3.57) from (3.55).
 From (3.56) and (3.57) we get
$$\ds\frac{H'(-\gamma_{j+1};j)}{H(\gamma_{j+1};j)}=
\ds\frac{H'(-\gamma_{j+1};0)}{H(\gamma_{j+1};0)}\ds\prod_{s=1}^j\left( \ds\frac{\gamma_{j+1}+\gamma_s}{\gamma_{j+1}-\gamma_s}
\right)^2
.\tag 3.58$$
Furthermore, from (3.36) and the fact that $F(k;0)$ has no zeros on the
positive imaginary axis, we know that $H(\beta;0)>0$ for $\beta>0.$
Thus, we see that (3.54) and (3.58) imply that
(3.38) and (3.53) are equivalent. \qed

One important consequence of Theorem~3.10 is that an ineligible resonance
remains ineligible if a number of bound states are removed or added via
Darboux transformations without changing the compact support of the potential.
An examination of the graph of $H(\beta;j)$ or $H(\beta;0)$ and the use of (3.38) or (3.53)
reveal various facts about eligible and ineligible resonances. The following proposition lists several such facts. We remind the reader that the meaning of
the maximal number of
eligible resonances is given in Section~1.

\noindent {\bf Proposition 3.11} {\it Let $V(\cdot; N)\in\Cal A$ be the potential of
the Schr\"odinger operator specified in (2.1)-(2.3)
with the boundary parameter $\theta$ in (2.3) and
$N$ bound states
at $k=i\gamma_j$ for $j=1,\dots,N,$ where we have $N=0$ if there are no bound states.
Let $M$ and $N_{\text{inel}}$ denote the maximal number of eligible resonances
and the number of ineligible resonances, respectively,
corresponding to the set $\{V(\cdot;N),\theta\}.$
Let $H(\beta;N)$ be the quantity corresponding to
the set $\{V(\cdot;N),\theta\},$ as defined
in (3.34). We have the following:}
%XXXXXXXXXXXXXXXX
%
%
%
%
%, i.e.
%the number of negative $\beta$-values satisfying
%$H(\beta;0)=0$ and $H'(\beta;0)>0,$ where i.e.
%the number of negative $\beta$-values satisfying
%$H(-\beta;0)=0$ and $H'(-\beta;0)\le 0.$
%}

\item{(a)} {\it The maximal
number of eligible resonances corresponding to the
set $\{V(\cdot;N),\theta\}$ is equal to the sum of the
number eligible resonances and the number of bound states
for $\{V(\cdot;N),\theta\}.$}

\item{(b)} {\it The number of ineligible resonances
for $\{V(\cdot;N),\theta\},$ i.e. the value of $N_{\text{inel}},$
remains unchanged if any number of bound states are removed or added
via Darboux transformations without changing the compact
support of the potential.}

\item{(c)} {\it Between any two consecutive eligible resonances
corresponding to $\{V(\cdot;N),\theta\},$
there must at least be
one ineligible resonance.}

\item{(d)} {\it We must have $M\le 1+N_{\text{inel}},$ and hence
for $\{V(\cdot;N),\theta\}$
we must also have $N\le 1+N_{\text{inel}}.$}

\item{(e)} {\it If there are at least two bound states associated with
the set $\{V(\cdot;N),\theta\},$ then there must at least be one ineligible resonance.}

\item{(f)} {\it
If $k=-i\gamma$ corresponds to an imaginary resonance and if $H(\beta;N)$
has no zeros in the interval $\beta\in(-\gamma,\gamma),$ then
$k=-i\gamma$ must correspond to an eligible
resonance for the set $\{V(\cdot;N),\theta\}.$}

\noindent PROOF: The proof of (a) intuitively follows from the definition of
the maximal number of eligible resonances, which is given in Section~1.
Here we provide the technical details.
Because $V(\cdot;N)\in \Cal A,$ by Theorem~2.1(a) the
corresponding Jost function $F(k;N)$ is entire in
$k\in\bC$ and hence $H(\beta;N)$ appearing
in (3.34) is a real-valued analytic function of $\beta\in\bR.$
By Theorems~2.1 and 3.3 it then follows that $H(\beta;s)$
is also a real-valued analytic function of $\beta\in\bR$ for any $s=0,1,\dots,N.$
By definition, $H(\beta;s)$ has exactly $s$ zeros in the interval
$\beta\in(0,+\infty),$ and by (3.53)
we conclude that $M$ is the number of zeros of $H(\beta;0)$
in the interval $\beta\in(-\infty,0)$ satisfying
$H'(\beta;0)>0.$ Thus, $H(\beta;N)$ is obtained from
$H(\beta;0)$ by converting $N$ eligible resonances into
bound states.
Hence, $H(\beta;N)$ has exactly $N$ bound states and $M-N$
eligible resonances, proving (a). From (3.53) it follows
that an ineligible resonance for $\{V(\cdot;N),\theta\}$
corresponds to a zero of the associated
 $H(\beta;0)$ in the interval $\beta\in(-\infty,0)$ satisfying
$H'(\beta;0)\le 0.$ As bound states are added, no such zeros of $H(\beta;0)$
are moved from the interval $(-\infty,0)$ to the interval $(0,+\infty).$
 Hence, (b) holds.
Let us now consider (c) when there are no bound states
so that we can use the eligibility criteria (3.53)
of Theorem~3.10.
In that case,
$H(\beta;0)$
is a real-valued analytic function of $\beta$ in the interval $(-\infty,0),$
and hence it is impossible to have two consecutive zeros of $H(\beta;0)$
in the interval $(-\infty,0)$
at which $H'(\beta;0)>0.$ Thus, in the absence of bound states
there has to be at least one ineligible
resonance between two eligible resonances. As stated
in the proof of (b), the ineligible resonances are unaffected
 if some eligible resonances are converted into bound states. Thus, the process of
 adding bound states does not change the location of the ineligible
 resonances but only moves a number of eligible resonances into
 bound states. Hence, even in the presence of bound states,
 we must have at least one ineligible resonance between
 two consecutive eligible resonances, proving (c).
 Note that the first inequality in (d) directly follows from (c).
 By (a) we have $N\le M$ and hence the second
 inequality in (d) is a consequence of the first inequality in (d).
 Note that (e) directly follows from the second inequality in (d)
if we have $N\ge 2.$
We prove (f) as follows. If $H(-\gamma;N)=0$ we cannot have $H(\gamma;N)=0$
because otherwise the corresponding regular solution
$\varphi_\theta(k,x)$ given in (2.7) would have
to be identically zero at $k=i\gamma,$ contradicting (2.5).
Furthermore, if $H(-\gamma;N)=0$
and $H(\beta;N)$ has no zeros
in the interval $\beta\in(-\gamma,\gamma],$ then $H'(-\gamma;N)$
and $H(\gamma;N)$ must have the same sign. Hence, (3.38)
implies that $k=-i\gamma$ is an eligible resonance. \qed

\vskip 10 pt
\noindent {\bf 4. RECOVERY FROM THE SCATTERING MATRIX}
\vskip 3 pt

In this section we assume that we are given a scattering matrix $S_\theta(k)$ for $k\in\bR$ and
we know that $S_\theta$
comes from a
potential $V$ in class $\Cal A$ and from a boundary parameter $\theta$ for some $\theta\in(0,\pi],$ where $\theta$ appears in (2.3).
However, we do not know what $V$ is and we do not know what the value
of $\theta$ is. In fact we do not even know whether $\theta=\pi$ or $\theta\in (0,\pi).$ In other words, we are only given the continuous
part of the Marchenko data specified in (2.15)
and we only know the existence of $V$ in $\Cal A$ and the existence of $\theta\in(0,\pi].$ In this section we have two main goals. Our first main goal is to determine whether $S_\theta(k)$ uniquely determines both $V$ and $\theta.$ Our second main goal is to reconstruct $V$ and $\theta$ in the case of uniqueness, or to reconstruct all possible sets $\{ V,\theta\}$ yielding the
same scattering matrix $S_\theta$ in the case of nonuniqueness.

To help the reader to understand better the theory developed in this section,
we first summarize our findings:

\item{(i)} If the extension of $S_\theta(k)$ from $k\in\bR$ to $k\in\bC$ has at least one pole on the positive
imaginary axis, then $S_\theta$ uniquely determines $V$ and
$\theta.$ We present an explicit algorithm to reconstruct the corresponding
$V$ and $\theta$ from $S_\theta.$

\item{(ii)} If the extension of $S_\theta(k)$ from $k\in\bR$ to $k\in\bC$ has no poles
on the positive imaginary axis and we have $S_\theta(0)=-1,$ then $S_\theta$ uniquely determines $V$ and $\theta.$ We present an explicit algorithm to reconstruct the corresponding
$V$ and $\theta$ from $S_\theta.$

\item{(iii)} If the extension of $S_\theta(k)$ from $k\in\bR$ to $k\in\bC$ has no poles
on the positive imaginary axis and we have $S_\theta(0)=+1,$ then there are precisely
two distinct sets $\{V_1,\theta_1\}$ and $\{V_2,\theta_2\}$ corresponding to the
same $S_\theta.$ We have $\theta_1=\pi$ and $\theta_2=\pi/2,$
and the potentials $V_1$ and $V_2$ correspond to
 some full-line reflection coefficients $R(k)$ and $-R(k),$ respectively.
Neither of the two corresponding full-line Schr\"odinger operators
have any bound
 states, and they are both exceptional in the sense that
 $R(0)\ne -1.$ We present an algorithm to reconstruct
the sets $\{V_1,\theta_1\}$ and $\{V_2,\theta_2\}.$

We already know from Theorem~2.1(f) that $S_\theta(0)$ must be either
$-1$ or $+1.$ Thus, the three cases listed above cover all
possible scenarios.
Having summarized our findings we now present the theory yielding the results in (i), (ii), and (iii), starting with case (i).

\noindent{\bf {Case (i)}} Given $S_\theta(k)$ for $k\in\bR,$ by the uniqueness
of the meromorphic extension, the poles of $S_\theta(k)$ on the positive imaginary
axis are uniquely determined. We already know from Theorem~2.1(e)
that such poles must be simple.
Let us assume that there are $N$ such poles and they occur at $k=i\gamma_s$
for $s=1,\dots,N.$ For the unique reconstruction of $V$ and
$\theta,$ we proceed as follows:

\item{(a)} We record the set $\{\gamma_1,\dots,\gamma_N\}$ as input to the
Marchenko method in (2.18)-(2.20) toward the identification of the bound states.

\item{(b)} Next, we evaluate the residues
$\text{Res}(S_\theta,i\gamma_s)$ for $s=1,\dots,N;$ i.e., we uniquely
determine the residue of $S_\theta(k)$ at each bound-state pole at $k=i\gamma_s.$ We then look at the sign of
$i\,\text{Res}(S_\theta,i\gamma_s)$ for any one value of $s.$ With the help of (2.26),
if that sign is positive then we conclude that $\theta=\pi,$ and
if that sign is negative then we conclude that $\theta\in(0,\pi).$

\item{(c)} From the previous step we know whether we have $\theta=\pi$ or
$\theta\in (0,\pi).$ Then, we use the appropriate line
in (2.18) and
the corresponding set $\{S_\theta,\{\gamma_s,m_s\}_{j=s}^N\}$ in the Marchenko
procedure outlined in Section~2 and we uniquely determine $V$ as in (2.20).
In case $\theta\in (0,\pi),$ we use (2.24) to determine the value of $\theta.$

\noindent{\bf {Case (ii)}} Given $S_\theta(k)$ for $k\in\bR$ with $S_\theta(0)=-1$ and with the further knowledge that the extension of $S_\theta(k)$ from $k\in\bR$ to $k\in\bC$
does not have any poles on the positive imaginary axis, we proceed as follows. From the Marchenko theory
outlined in (2.18)-(2.24), we see that we only need to know whether we have
$\theta=\pi$ or $\theta\in (0,\pi).$ This is because we will use
either the first line or the second line of (2.18), but without the summation
terms in those lines, as input to the corresponding Marchenko equation.
Thus, in the Marchenko equation (2.19) we have the Marchenko kernel and the nonhomogeneous term are determined up to a sign, depending on whether we have $\theta=\pi$ or $\theta\in(0,\pi).$
Let us assume that
corresponding to $S_\theta,$ we have two distinct sets $\{V_1,\theta_1\}$ and $\{V_2,\theta_2\}.$ We cannot have both $\theta_1$ and $\theta_2$ equal to
$\pi$ because then the second line of
(2.18) would yield $V_1\equiv V_2$ via the Marchenko method. Similarly, we cannot have
both $\theta_1$ and $\theta_2$ different from $\pi$ because then the first
line of (2.18) would yield $V_1\equiv V_2.$ Thus, one of $\theta_1$ and $\theta_2$
must be equal to $\pi$ and the other must be different from $\pi.$ Without
loss of any generality we can assume that $\theta_1=\pi$ and $\theta_2\in (0,\pi).$
Let us use $f_1(k,0)$ to denote the Jost function corresponding to
$\{V_1,\theta_1\}$ and use $F_2(k)$ to denote the Jost function corresponding to
$\{V_2,\theta_2\}.$ Because $S_\theta(0)=-1,$ from (2.27) and (2.28) it follows that
we must have
$f_1(0,0)=0$ and $F_2(0)\ne 0.$ From Theorem~2.1 we know that $k=0$ must be a simple
zero of $f_1(k,0)$ and hence we have $f_1(k,0)=k\,h_1(k)$ for some function
$h_1(k)$ in such a way that $h_1(k)$ is analytic and nonzero in $k\in\bCpb$ and
$h_1(k)=1/k+O(1/k^2)$ as $k\to\infty$ in $k\in\bCpb.$
Similarly, from Theorem~2.1 we know that
$F_2(k)$ is analytic and nonzero in $k\in\bCpb$ and
$F_2(k)=k+O(1)$ as $k\to\infty$ in $k\in\bCpb.$
Since $f_1(k,0)$ and $F_2(k)$ correspond to the
same scattering matrix $S_\theta(k),$ because of (2.16) we must have
$$S_\theta(k)=\ds\frac{f_1(-k,0)}{f_1(k,0)}=\ds\frac{-F_2(-k)}{F_2(k)},
\qquad k\in\bR,\tag 4.1$$
which implies
$$\ds\frac{f_1(k,0)}{F_2(k)}=\ds\frac{-f_1(-k,0)}{F_2(-k)},\qquad k\in\bR.\tag 4.2$$
Since $f_1(k,0)=k\,h_1(k),$ we can write (4.2) also as
$$\ds\frac{h_1(k)}{F_2(k)}=\ds\frac{h_1(-k)}{F_2(-k)},\qquad k\in\bR.\tag 4.3$$
Note that the left-hand side of (4.3) has
 an analytic extension from $k\in\bR$ to
$k\in\bCp,$ and that analytic extension is
continuous in $\bCpb$ and behaves as $O(1/k^2)$ as $k\to \infty$
in $\bCpb.$ Similarly, the right-hand side of (4.3) has an analytic extension from $k\in\bR$ to
$k\in\bCm,$ and that analytic extension is continuous in $\bCmb$ and
behaves as $O(1/k^2)$ as $k\to \infty$
in $\bCmb.$ Thus, $h_1(k)/F_2(k)$ must be an entire function of $k$ and behaving
like $O(1/k^2)$ as $k\to \infty$
in $\bC.$ By Liouville's theorem, we must then have $h_1(k)\equiv 0.$ However,
that would imply $f_1(k,0)\equiv 0,$ contradicting the second line of (2.27).
Thus, we cannot have both $\{V_1,\theta_1\}$ and $\{V_2,\theta_2\}$ corresponding
to the same $S_\theta(k)$ and we must have a unique set $\{V,\theta\}$ corresponding to $S.$
Having established the uniqueness, let us now consider the reconstruction problem. As explained in Section~2, we can use the Marchenko
method for the reconstruction. We can first try the second line
of (2.18) as input to the Marchenko equation with
$\theta=\pi$ without the summation term there.
We can construct the corresponding
potential and Jost solution via (2.20) and can check
if the right-hand side of (2.21) is zero at $k=0,$
which is required by the second line of (2.27).
Alternatively, we can check if the right-hand side of (2.23)
is equal to our scattering matrix $S_\theta(k).$
If there is no agreement, we then know that $\theta\in (0,\pi),$ and
hence use the first line of (2.18) without the summation term there
as input to the Marchenko equation
and uniquely construct the corresponding $V$ and $\theta$ via the first
equality in (2.20) and by using (2.24), respectively.

\noindent{\bf {Case (iii)}} Given $S_\theta(k)$ for $k\in\bR$ with $S_\theta(0)=+1$ and with the further knowledge that the extension of $S_\theta(k)$ from $k\in\bR$ to $k\in\bC$
does not have any poles on the positive imaginary axis, we proceed as follows.
As in case (ii), from the Marchenko theory it follows that it is enough to check the
nonuniqueness by assuming that,
corresponding to $S_\theta,$ we have two distinct sets $\{V_1,\theta_1\}$ and $\{V_2,\theta_2\}$ with $\theta_1=\pi$ and $\theta_2\in (0,\pi).$
Contrary to case (ii), we will now prove that there are precisely two distinct sets
$\{V_1,\theta_1\}$ and $\{V_2,\theta_2\}$ corresponding to the same $S_\theta(k).$
We again use $f_1(k,0)$ to denote the Jost function corresponding to
$\{V_1,\theta_1\}$ and use $F_2(k)$ to denote the Jost function corresponding to
$\{V_2,\theta_2\}.$
Let us use $f_2(k,x)$ to denote the Jost solution
corresponding to $V_2.$ From (2.6) we
have
$$F_2(k)=-i\left[f_2'(k,0)+(\cot\theta_2)\, f_2(k,0)\right].\tag 4.4$$
This time, from (2.27) and (2.28) it follows that
$f_1(0,0)\ne 0$ and $F_2(0)=0.$ From Theorem~2.1(a) we know that
$k=0$ must be a simple
zero of $F_2(k),$ and hence
we have $F_2(k)=k\,g_2(k)$ for some function
$g_2(k)$ in such a way that $g_2(k)$ is analytic and nonzero in $k\in\bCpb$ and
$g_2(k)=1+O(1/k)$ as $k\to\infty$ in $k\in\bCpb.$
Similarly, from Theorem~2.1 we know that
$f_1(k,0)$ is analytic and nonzero in $k\in\bCpb$ and
$f_1(k,0)=1+O(1/k)$ as $k\to\infty$ in $k\in\bCpb.$
Since $f_1(k,0)$ and $F_2(k)$ correspond to the
same scattering matrix $S_\theta(k),$ we must have (4.1) and (4.2) satisfied.
Since $F_2(k)=k\,g_2(k),$ we can write (4.2) also as
$$\ds\frac{f_1(k,0)}{g_2(k)}=\ds\frac{f_1(-k,0)}{g_2(-k)},
\qquad k\in\bR.\tag 4.5$$
Note that the left-hand side of (4.5) has an analytic extension from $k\in\bR$ to
$k\in\bCp,$ and that analytic extension
is continuous in $\bCpb$ and behaves as $1+O(1/k)$ as $k\to \infty$
in $\bCpb.$ Similarly, the right-hand side of (4.5) has an analytic extension from $k\in\bR$ to
$k\in\bCm,$ and that analytic extension is continuous
in $\bCmb$ and behaves like $1+O(1/k)$ as $k\to \infty$
in $\bCmb.$ Thus, we must have $f_1(k,0)/g_2(k)$ entire and behaving
like $1+O(1/k)$ as $k\to \infty$
in $\bC.$ By Liouville's theorem, we must then have $g_2(k)\equiv f_1(k,0),$
or equivalently we must have
$$F_2(k)\equiv k\,f_1(k,0).\tag 4.6$$
Since there are no poles of $S_\theta(k)$ on the positive imaginary axis in $\bCp,$
it follows that the Marchenko kernel, which we call $M_1(y),$
corresponding to the first set $\{V_1,\theta_1\}$
with $\theta_1=\pi$ is given by the second line of (2.18) but without the summation
term there. Then, the Marchenko kernel, which we call $M_2(y),$ corresponding to the second set $\{V_2,\theta_2\}$
with $\theta_2\in (0,\pi)$ is given by the first line of (2.18) but without the summation
term there. From (2.18) it is clear that $M_2(y)=-M_1(y).$ Let us now view
$V_1$ and $V_2$ as compactly-supported potentials in the full-line Schr\"odinger
equation with $V_1(x)\equiv 0$ for $x<0$ and $V_2(x)\equiv 0$ for $x<0.$
As in (2.30) let us associate the scattering
coefficients $T_1,$ $L_1,$ $R_1$ with
$V_1$ and associate the scattering coefficients $T_2,$ $L_2,$ $R_2$ with
$V_2.$
Since $S_\theta(k)$ has no poles on the positive imaginary axis, we know that
$N_{\theta_2}=0$ and $N_\pi=0,$ where $N_{\theta_2}$ and
$N_\pi$ denote the number of bound states corresponding to
$\{V_2,\theta_2\}$ and $\{V_1,\theta_1\},$ respectively.
 From (4.6) we know that $F_2(0)=0,$ and hence Theorem~2.2(f)
indicates that we cannot have $\theta_2\in (0,\pi/2)\cup (\pi/2,\pi]$ and thus
we must have $\theta_2=\pi/2,$ which yields $\cot\theta_2=0.$
By Theorem~2.3 we then know that $\tilde N=0,$ i.e. neither $T_1(k)$ nor $T_2(k)$
has any poles
on the positive imaginary axis. From Theorem~2.4 we then get $M_1(y)=\hat R_1(y)$ and
$M_2(y)=\hat R_2(y)$ for
$y>0$ with $\hat R_1(y)$ and $\hat R_2(y)$ denoting the Fourier
transforms as in (2.45). Since we already know that $M_2(y)\equiv -M_1(y),$
we then get $\hat R_2(y)\equiv -\hat R_1(y),$ and hence yielding
$R_2(k)\equiv -R_1(k).$ Because $\tilde N=0,$ it then also follows that
$T_1(k)\equiv T_2(k).$ From the characterization conditions [8-10,14,15] for the full-line
Schr\"odinger operators, we already know that if there exists $V_1\in \Cal A$ corresponding to
$R_1$ and $T_1,$ we are assured the existence of $V_2\in\Cal A$ corresponding
to $-R_1$ and $T_1$ by recalling that $R_1(0)\ne -1$ and that $T_1$ does not have any
bound-state poles on the positive imaginary axis. Thus, we have established
the existence of two distinct sets $\{V_1,\theta_1\}$ and $\{V_2,\theta_2\}$
with $\theta_1=\pi$ and $\theta_2=\pi/2.$
Note that, using $\cot\theta_2=0$ in (4.4) we get $F_2(k)=-i f_2'(k,0)$
and hence (4.6) indicates that
$$f_2'(k,0)=ik\,f_1(k,0).$$
One consequence of (4.6) is that we must have
$$\int_0^b dx\,V_1(x)=\int_0^b dx\,V_2(x).\tag 4.7$$
We obtain (4.7) by expanding $F_2(k)$ with $\cot\theta_2=0$
with the help of the first line of (2.25) and by comparing it with
 the expansion of the right-hand side of (4.6) via the second line
 of (2.25).
The potential $V_1$ can be reconstructed with the help
 of the second line of (2.18) without the summation term there.
 The potential $V_1$ is then obtained by solving (2.19) and
using the first equality in (2.20). Similarly, $V_2$ can be reconstructed by using the first line of (2.18)
without the summation term there. Thus, $V_2$ can be
obtained by using (2.19) and (2.20).
%The scattering coefficients $R_1,$ $L_1,$ $T_1$ have various other restrictions on them
%if $-R_1,$ $-L_1,$ $T_1$ corresponds to the same scattering matrix but
%coming from another
% half-line Schr\"odinger
%operator. For example, from the first identity in (2.30) and the first equality in
%(4.1) we obtain
%$$\ds\frac{\ds\frac{1+L_1(-k)}{T_1(-k)}}{\ds\frac{1+L_1(k)}{T_1(k)}}=
%\ds\frac{\ds\frac{1-L_1(-k)}{T_1(-k)}}{\ds\frac{1-L_1(k)}{T_1(k)}},\tag 4.8$$
%which simplifies to
%$$L_1(-k)=L_1(k),\qquad k\in\bC.\tag 4.9$$
%For $k\in\bR,$ it is already known that $L_1(-k)=L_1(k)^\ast,$ where the asterisk
%denotes complex conjugation. Thus, $L_1(k)$ must be real valued for $k\in\bR.$
%Using the third identity of (4.8) in (4.9) we get
%$$R_1(-k)=R_1(k)\,\ds\frac{T_1(-k)^2}{T_1(k)^2},\qquad k\in\bC.$$

We summarize our findings in this section in the following theorem.

\noindent {\bf Theorem 4.1} {\it Assume that we are given $S_\theta(k)$ for
$k\in\bR$ and we know that it comes from a
potential $V$ in class $\Cal A$ and from a boundary parameter $\theta$ for some $\theta\in(0,\pi],$ where $\theta$ appears in (2.3).
We then have the following:}

\item{(a)} {\it If $S_\theta(0)=+1$ and
the extension of $S_\theta(k)$ from $k\in\bR$
to $k\in\bCp$ has no poles on the
positive imaginary axis, then there are precisely
two distinct sets $\{V_1,\theta_1\}$ and $\{V_2,\theta_2\}$
with $\theta_1=\pi,$ $\theta_2=\pi/2,$ $V_1\in\Cal A,$ and $V_2\in\Cal A.$
The set $\{V_1,\theta_1\}$ corresponds to the Jost solution $f_1(k,x),$ and the
corresponding Jost function $f_1(k,0)$ satisfies $f_1(0,0)\ne 0.$ The Jost
function $F_2(k)$ for the second set $\{V_2,\theta_2\}$ is
equal to $k f_1(k,0).$
Both sets can be uniquely reconstructed by the Marchenko procedure. The
set $\{V_1,\theta_1\}$ is associated with some
scattering coefficients $R_1,$ $L_1,$ $T_1$ in such a way that
%$L_1(k)$ is an even function
%of $k\in\bC,$ $L_1(k)$ is real valued for
%$k\in\bR,$
$T_1(0)\ne 0$ and that $T_1(k)$ does not have poles on the
positive imaginary axis. The scattering coefficients $R_1,$ $L_1,$ $T_1$
are related to $f_1(k,0)$ and $f_1'(k,0)$ as in (2.30).
The set
$\{V_2,\theta_2\}$ is associated with the
scattering coefficients $R_2,$ $L_2,$ $T_2$
where $R_2(k)\equiv -R_1(k),$ $L_2(k)\equiv -L_1(k),$ and $T_2(k)\equiv T_1(k).$
Although, in general the potentials $V_1$ and $V_2$
are distinct, their
integrals have the same value, as seen from (4.7). The very special case
$V_1(x)\equiv V_2(x)$ occurs when $R_1(k)\equiv 0,$ $L_1(k)\equiv 0,$ $T_1(k)\equiv 1,$
which yields $V_1(x)\equiv 0$ and $V_2(x)\equiv 0$.}

\item{(b)} {\it If $S_\theta(0)\ne +1$ or
the extension of $S_\theta(k)$ from $k\in\bR$
to $k\in\bCp$ has at least one pole
on the
positive imaginary axis, then there is a unique
potential $V\in\Cal A$ and a unique boundary parameter
$\theta$ in the interval $(0,\pi]$ corresponding
to $S_\theta(k).$ The corresponding potential $V$ and boundary parameter $\theta$
can be uniquely reconstructed by the Marchenko procedure
outlined in Section~2.}

\vskip 10 pt
\noindent {\bf 5. RECOVERY FROM ABSOLUTE VALUE OF THE JOST FUNCTION}
\vskip 3 pt

Our goal in this section is to investigate the determination of a real-valued,
integrable, compactly-supported
potential and a selfadjoint boundary condition from the input data
consisting of the absolute value of the corresponding Jost function
known at positive energies.
In other words,
we assume that we only know the continuous part of the Gel'fand-Levitan spectral data
given in
(2.10) without having any explicit knowledge of its discrete part. Furthermore, we know that our input data set corresponds to a selfadjoint Schr\"odinger operator on the half line with a selfadjoint
boundary condition at $x=0.$ However, we do not know if the boundary condition
is Dirichlet or non-Dirichlet, and
we do not know if there are any bound states and we do not know the number of bound states
if there are any. In fact, we would
like to determine all such characteristics
 from our input data set alone, if possible.

In this section we use the notation introduced in Section~3, namely,
we use
$\theta_j,$ $V(x;j),$ $\varphi(k,x;j),$ and
$F(k;j)$ to denote the relevant
quantities corresponding to the half-line Schr\"odinger operator
with bound states at $k=i\gamma_1,\dots,i\gamma_j,$ where
the case $j=0$ refers to the quantities
with no bound states. Note that $\theta_j$ is the boundary
parameter appearing in (2.3), $\varphi(k,x;j)$ is the regular solution
in (2.5), $F(k;j)$ is the Jost function in (2.6),
$g_j$ is the Gel'fand-Levitan norming constant in (2.9), $G(x,y;j)$ is the
Gel'fand-Levitan kernel appearing in (2.12) and (2.13),
$A(x,y;j)$ is the solution in (2.14) to the Gel'fand-Levitan equation,
and $H(\beta;j)$ is the quantity in (3.34).

Mathematically speaking, we consider the selfadjoint Schr\"odinger operator on the half line
with the
potential $V(\cdot;N)$ in class $\Cal A,$ the boundary parameter $\theta_N,$
the Jost function $F(k;N),$ the bound states at $k=i\gamma_s$ with the
corresponding Gel'fand-Levitan norming constants $g_s$ for $s=1,\dots,N,$
where $N$ is a nonnegative integer. We assume that our input data set solely
consists of $|F(k;N)|$ for $k\in\bR.$ We do not know
the value of $N,$ and we do not
know anything about the set $\{\gamma_s,g_s\}_{s=1}^N.$ We would like to investigate
to what extent our input data set determines $N,$ $\theta_N,$ $\{\gamma_s,g_s\}_{s=1}^N,$ and $V(x;N).$ In other words, we know the existence
of at least one potential $V$ in class $\Cal A$ and the existence
of one selfadjoint boundary parameter
$\theta\in (0,\pi]$ corresponding to our input data, and we
would like to investigate the uniqueness or
nonuniqueness of the set $\{V,\theta\}$
by determining all potentials $V$
in class $\Cal A$ and all boundary parameters $\theta$ in the interval
$(0,\pi]$
 corresponding to our
input data set.

Our findings are summarized as follows:
We can uniquely determine whether the boundary condition
is Dirichlet or non-Dirichlet. We can determine all the
corresponding potentials and
boundary conditions, but the uniqueness is only up to the inclusion of the
eligible resonances.
Thus, if the maximal number of eligible resonances is zero, then we have
the unique determination of the potential $V(x;0)$ and the boundary
parameter $\theta_0$ corresponding to our data. If the maximal number of
eligible resonances is one, then we determine the two distinct sets
$\{V(x;0),\theta_0\}$ and $\{V(x;1),\theta_1,\{\gamma_1,g_1\}\}$ corresponding
to our input data. If the maximal number of eligible
resonances is $M,$ then we determine that there is a $2^M$-fold nonuniqueness
and that any one of those $2^M$ sets corresponds to our input data.
We remind the reader that the definition of
the maximal number of eligible resonances is given in Section~1.

As mentioned earlier, the number of imaginary resonances may be infinite, but
under some mild additional assumptions [19] such as $V(x)\ge 0$ or $V(x)\le 0$ in the vicinity
of $x=b,$ that number is guaranteed to be finite.
We recall that $b$ refers to the constant in (2.2)
and related to the compact support of the potential
$V.$ Thus, under a mild
additional assumption we are guaranteed that $M,$ the maximal
number of eligible resonances, is finite.

Having summarized our findings, let us now outline the method of
determining all potentials and boundary conditions corresponding to our input data:

\item{(a)} From our input data $|F(k;N)|$ for $k\in\bR,$ by using
the asymptotic behavior in (2.11) we can tell whether
the corresponding boundary parameter $\theta_N$ satisfies
$\theta_N\in (0,\pi)$ or
$\theta_N=\pi.$

\item{(b)} From (3.3) and (3.6) it is clear that we have
$$|F(k;0)|=|F(k;N)|,\qquad k\in\bR,\tag 5.1$$
where $F(k;0)$ is the Jost function corresponding to no bound states. Using the
Gel'fand-Levitan procedure outlined in Section~2, from $|F(k;0)|,$ which
is equivalent to $|F(k;N)|$ as seen from (5.1), we uniquely
construct $V(x;0),$ $\theta_0,$ and the regular solution $\varphi(k,x;0).$
This is done, by first forming the Gel'fand-Levitan kernel as in (2.12)
and (2.13), namely
$$G(x,y;0):=\cases \ds\frac{1}{\pi}\ds\int_{-\infty}^\infty dk\,\left[
\ds\frac{k^2}{|F(k;N)|^2}-1\right] (\cos kx)(\cos ky),\qquad & \theta_N\in (0,\pi),\\
\stretch
\ds\frac{1}{\pi}\ds\int_{-\infty}^\infty dk\,\left[
\ds\frac{1}{|F(k;N)|^2}-1\right] (\sin kx)(\sin ky),\qquad & \theta_N=\pi.\endcases$$
Using $G(x,y;0)$ in the corresponding Gel'fand-Levitan equation in (2.14), namely
in
$$A(x,y;0)+ G(x,y;0)+\ds\int_0^x dz\, A(x,z;0)\,G(z,y;0)=0,
\qquad 0<y<x,$$
we uniquely recover $A(x,y;0),$ from which we get $V(x;0),$ $\theta_0,$ and $\varphi(k,x;0)$ via
$$\cot \theta_0=-A(0,0;0),\qquad \theta_N\in(0,\pi),$$
$$V(x;0)=2\,\ds\frac{d A(x,x;0)}{dx},\qquad \theta_N\in(0,\pi].$$

\item{(c)} As a consequence of (5.1),
we uniquely determine $F(k;0)$ from our input data via [7]
$$F(k;0)=\cases
k\,\exp\left(
\ds\frac{-1}{\pi i}\ds\int_{-\infty}^\infty
dt\,\ds\frac{\log|t/F(t;N)|}{t-k-i0^+}
\right),&\qquad \theta_N\in(0,\pi),\\
\exp\left(
\ds\frac{1}{\pi i}\ds\int_{-\infty}^\infty
dt\,\ds\frac{\log|F(t;N)|}{t-k-i0^+}
\right),&\qquad \theta_N=\pi,\endcases\tag 5.2$$
where $i0^+$ indicates that the value for $k\in\bR$ must be obtained as a
limit from within $\bCp.$ Since $F(k;0)$ has an analytic
 extension to the entire complex plane,
we are assured that
(5.2) holds for all $k\in\bC.$

\item{(d)} Having $F(k;0)$ at hand for $k\in\bC,$ we construct the real-valued
function $H(\beta;0)$ defined in (3.34).
We already know that $H(\beta;0)$ does not have any zeros when $\beta>0.$
We can have $H(0;0)\ne 0$ (generic case) or we can have $H(0;0)=0$ (exceptional case)
with a simple zero of $H(\beta;0)$ at $\beta=0.$ We then go ahead and determine
all imaginary resonances, i.e. the zeros of $H(\beta;0)$ when $\beta<0.$

\item{(e)} We then identify each imaginary resonance either as eligible or
ineligible by using the eligibility criteria given in (3.53), namely
by finding all negative $\beta$-values satisfying
$$H(\beta;0)=0,\quad H'(\beta;0)> 0.\tag 5.3$$
Assuming that (5.3) is satisfied when $\beta=-\beta_s$ for $s=1,\dots,M,$
we uniquely determine the set $\{\beta_s\}_{s=1}^M.$ Note that $M$ is the maximal
number of eligible resonances. We know that $M$ may be zero, a positive
integer, or infinity. As mentioned previously, a mild additional assumption [19]
guarantees the finiteness of $M.$

\item{(f)} Each eligible resonance $k=-i\beta_s$ can be converted into a bound state
by using the Darboux transformation formulas given in Theorem~3.1. Thus, it is possible
to add $N$ bound states, where $N$ is an integer between $0$ and $M.$ We can choose
$N$ bound states at $k=i\gamma_s$ among the $M$ possible choices
$k=i\beta_s$ in $\binom{M}{N}$ ways, where $\binom{M}{N}$ denotes the binomial coefficient,
which is equal to $M!/((N!)(M-N)!).$ Thus, as $N$ takes all values
between $0$ and $M,$ we find that we have
$2^M$ distinct sets consisting of a potential and a boundary parameter, each corresponding to the
same absolute value of the Jost function.

\vskip 10 pt
\noindent {\bf 6. EXPLICIT EXAMPLES}
\vskip 3 pt

In this section we illustrate our main results
presented in Sections~3-5 with some explicit examples.
The first example is provided to remind the reader
that the boundary parameter $\theta$
 appearing in (2.3)
 indeed affects the bound states and resonances, and in fact even the
trivial potential can have a bound state or a resonance depending on the
value of the boundary parameter $\theta$
appearing in (2.3).

\noindent {\bf Example 6.1} Assume that $V(x)\equiv 0$ in (2.1). The corresponding Jost function $F_\theta(k)$
is given by (2.8).
Since $F_\pi(k)$ has no zeros in $\bC,$
there are no bound states and there are
no resonances in the Dirichlet case $\theta=\pi.$
Let us now consider the non-Dirichlet case with some fixed boundary parameter $\theta\in(0,\pi).$
Recall that the zeros of $F_\theta(k)$ in $\bCp$ correspond to the
bound states and the zeros in $\bCm$ correspond to the resonances.
If $\cot\theta>0,$ then there is one bound state and there are no resonances. If $\cot\theta=0,$
then there are no bound states and there are no resonances.
If $\cot\theta<0,$ then there are no bound states and there is one imaginary resonance.
In fact, as a result of Proposition~3.11(f), $k=i\cot\theta$ is an eligible resonance
when $\cot\theta<0.$
Thus,
if $\cot\theta<0$ we can add a bound state to $V\equiv 0$
at $k=-i\cot\theta,$ and if we choose the Gel'fand-Levitan
bound-state norming constant $g$ as in (3.19), i.e.
with $g^2=-2\cot\theta,$ then the transformed potential
still vanishes everywhere, and hence the transformed
potential and the original potential have the same (trivial)
compact support. Note that such a choice is compatible with (3.1). Let us see
what happens if we do not use $g^2=-2\cot\theta$ as our norming
constant. With $f(k,x)= e^{ikx}$ and $F_\theta(k)=k-i\,\cot\theta,$
using the first line in (2.7)
 we evaluate $\varphi_\theta(k,x)$ as
 $$\varphi_\theta(k,x)=\ds\frac{1}{2k}\left[(k-i\,\cot\theta)
 \,e^{-ikx}+(k+i\,\cot\theta)
 \,e^{ikx}\right].$$
If we add a bound state at $k=-i\,\cot\theta$ with the Gel'fand-Levitan
norming constant $g,$ then the quantity inside the brackets in (3.2) is given by
the right-hand side in the following equation:
$$\ds\frac{2 g^2\, \varphi_\theta(-i\,\cot\theta,x)^2}{
1+ g^2\,\ds\int_0^x dy\,\varphi_\theta(-i\,\cot\theta,y)^2}
=\ds\frac{4 g^2 \cot\theta}{-g^2+(2 \cot\theta+g^2)\,\ds e^{ 2 x\cot\theta}}.\tag 6.1$$
Thus, the choice $g^2=-2 \cot\theta$ makes the right-hand side in
(6.1) equal to the constant $-4\cot\theta,$ and hence the support
of the potential is unchanged when we add the bound state
at $k=-i\,\cot\theta$ with the norming constant $g=\sqrt{-2 \cot\theta}.$
Any other choice for the norming constant $g$
results in a potential with support on the entire half line.

Next, we provide some examples of
eligible resonances when the potential and the boundary parameter
are known.

\noindent {\bf Example 6.2} Let us assume that we are given the boundary parameter
$\theta\in(0,\pi)$ and that $V(x)$ is
the piecewise constant potential (potential barrier or potential well) given by
$$V(x)=\cases v,\qquad 0<x<1,\\
\stretch
0,\qquad x>1,\endcases\tag 6.2$$
where $v$ is a constant parameter.
With the help of (2.4)-(2.7) and (6.2) we
can explicitly evaluate the regular solution
$\varphi_\theta(k,x),$ the Jost solution $f(k,x),$ and the Jost
function $F_\theta(k)$ and get
$$\varphi_\theta(k,x)=\cases
\cosh \eta x-\cot\theta\,\ds\frac{\sinh \eta x}{\eta},
&\qquad 0\le x\le 1,\\
\stretch
b_1\,\cos k(x-1)+\ds\frac{b_2}{k}\,\sin k(x-1),&\qquad 1\le x< +\infty,
\endcases\tag 6.3$$
$$f(k,x)=
\cases
e^{ik}\,\cosh \eta (x-1)+ik\,e^{ik}
\,\ds\frac{\sinh \eta (x-1)}{\eta},
&\qquad 0\le x\le 1,\\
\stretch
e^{ikx},&\qquad 1\le x< +\infty,
\endcases$$
$$F_\theta(k)=e^{ik}
(k-i\cot\theta)\cosh\eta-e^{ik}(k\cot\theta-i\eta^2)\,\ds\frac{\sinh \eta}{\eta}
,\tag 6.4$$
where we have defined
$$\eta:=\sqrt{v-k^2},\quad b_1:=\cosh\eta-\cot\theta\,\ds\frac{\sinh \eta}{\eta},
\quad b_2:=\eta\, \sinh \eta-\cot\theta\,\cosh\eta.$$
Let us now analyze (6.2) for various values of $v$ and $\cot\theta.$
We use an overline on a digit
to indicate a round off.

\item{(a)} When $(v,\cot\theta)=(-10,1),$ using (3.34)
and (6.4) we obtain $H_\theta(\beta),$ plotted in the
first graph of Figure~6.1.
We
observe from the graph of $H_\theta(\beta)$
that it has two positive zeros and one negative zero. Thus,
there are two bound states occurring at
$k=0.76040\overline{9}i$ and $k=3.2527\overline{3}i$ and
that $F_\theta(k)$ has a simple zero at $k=-\gamma i,$
where $\gamma=2.8208\overline{4}.$
 From the graph of
$H_\theta(\beta)$ we easily see that
$H_\theta(\gamma)<0$ and $H'_\theta(-\gamma)>0,$
and hence by (3.38) we conclude
that $k=-\gamma i$ is an ineligible resonance and that it is impossible
to add a bound state to $V$ without changing the compact
support property. Equivalently, using $b=1$ for the constant $b$
appearing in (2.2), with the help of
(6.3) we evaluate the right-hand side of the
second equality in (3.19) and hence obtain
$g^2=-4.2376\overline{1}.$ Thus, we confirm that
$k=-\gamma i$ is an ineligible resonance because (3.20) is not satisfied.
The same conclusion can also be reached via
Proposition~3.11(e) because we have precisely two bound states and one imaginary resonance
and hence that imaginary resonance must be ineligible.

\item{(b)} When $(v,\cot\theta)=(-0.2,6),$ the
plot of $H_\theta(\beta),$ given as the second graph in
Figure~6.1, reveals that
$H_\theta(\beta)$ has one positive zero and two negative zeros. Thus,
there is a bound state at
$k=6.0166\overline{4}i$ and that
$F_\theta(k)$ has simple zeros at $k=-\gamma_1 i$
and $k=-\gamma_2 i,$
where $\gamma_1=3.3618\overline{2}$
and $\gamma_2=5.9584\overline{2}.$ From the graph of
$H_\theta(\beta)$ we easily see that
$H_\theta(\gamma_2)<0$ and $H'_\theta(-\gamma_2)>0,$
and hence $k=-\gamma_2i$ is an ineligible
resonance, as indicated by the criteria in (3.38).
 On the other hand,
$H_\theta(\gamma_1)<0$ and $H'_\theta(-\gamma_1)<0,$
so that $k=-\gamma_1i$ is an eligible resonance because of the criteria
in (3.38).
In fact from the second equality in (3.19), using $b=1$ and
$\gamma=\gamma_1$ we get
$g^2=g_1^2>0$ with $g_1^2=1.9320\overline{9}.$
Thus, we can add a bound state to $V$ at $k=i\gamma_1$
with the Gel'fand-Levitan norming constant $g_1=1.3\overline{9}$ and the
resulting potential has also support in the interval $(0,1).$

\item{(c)} When $(v,\cot\theta)=(0.003521,-3),$ from the
plot of $H_\theta(\beta)$ given as the third graph in Figure~6.1
we observe that
$H_\theta(\beta)$ has no positive zeros and has a double zero
at a negative $\beta$-value. Thus,
there are no bound states and
$F_\theta(k)$ has a double zero at $k=-\gamma i,$
where $\gamma=3.620\overline{5}.$
We have $H_\theta(\gamma)>0$ and $H'_\theta(-\gamma)=0.$
Thus, the incompatibility with (3.38) shows that
we cannot add any bound states to
$V$ without changing the compact support property.

\vskip 5 pt

\centerline{\hbox{
\psfig{figure=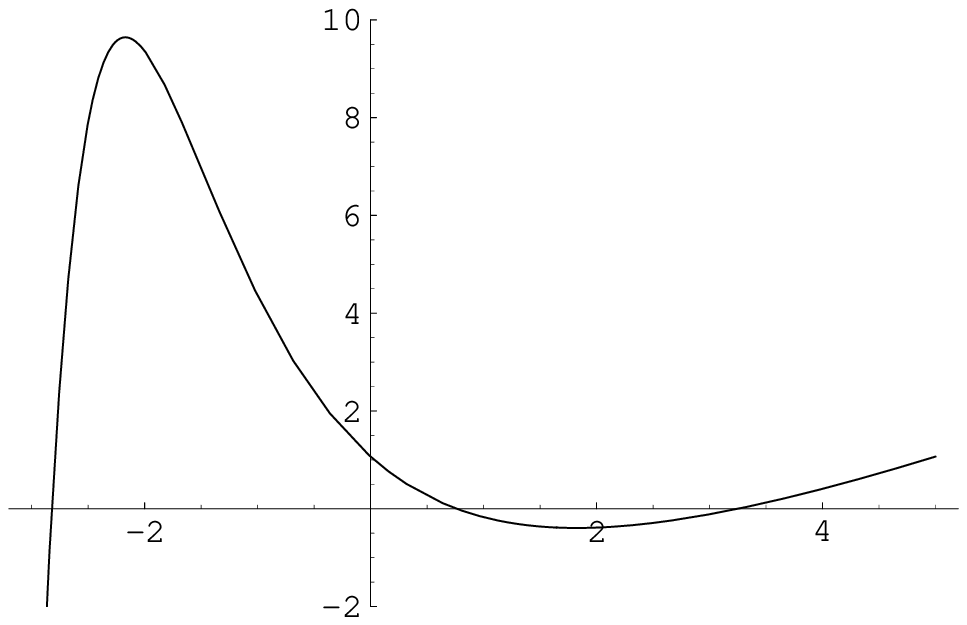,width=2 truein,height=2 truein}\quad
\psfig{figure=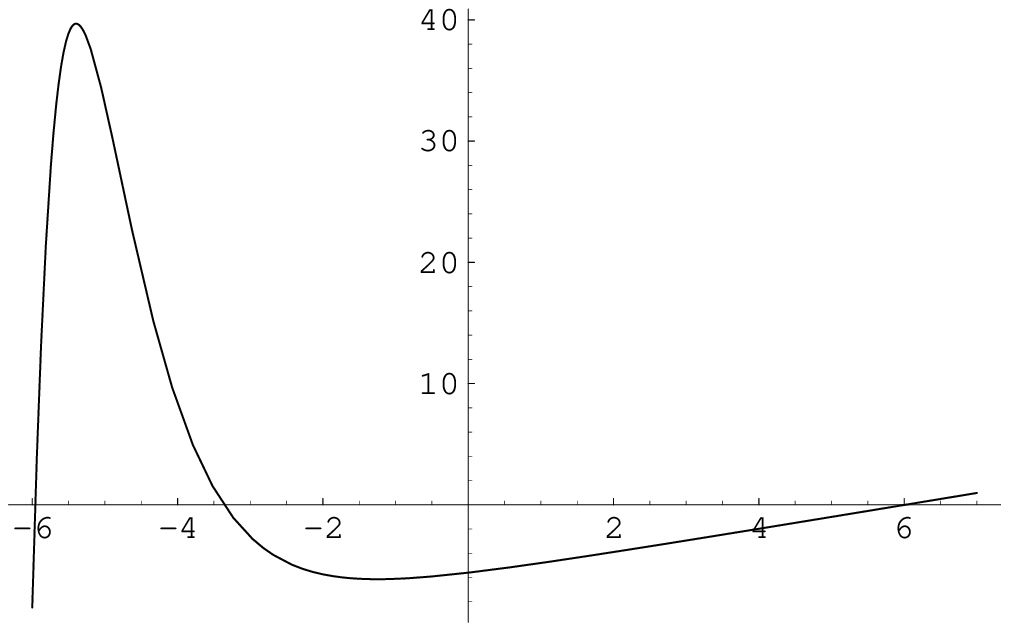,width=2.1 truein,height=2.1 truein}\quad
\psfig{figure=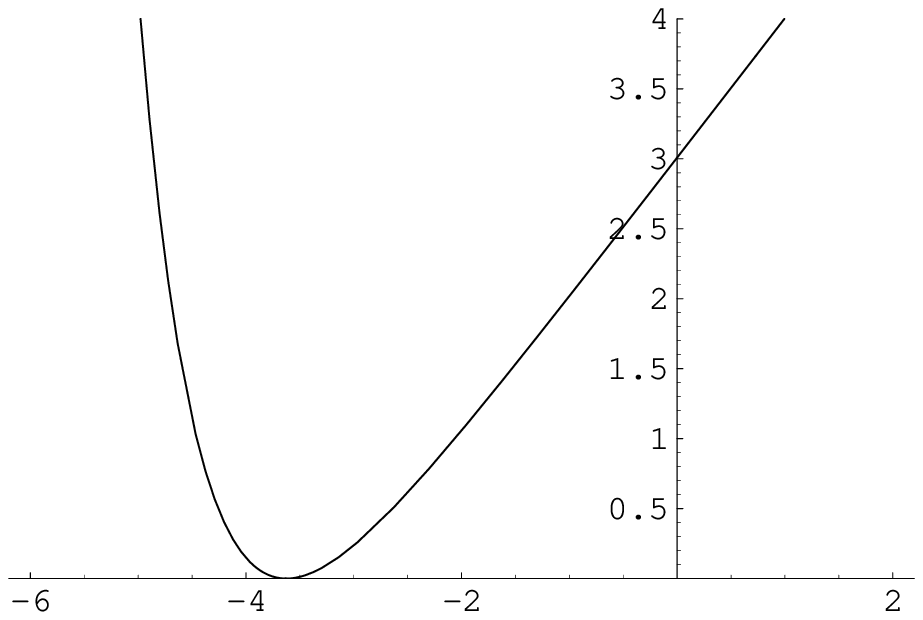,width=2.1 truein,height=2 truein}
}}

\centerline{{\bf Figure~6.1} The plots of $H_\theta(\beta)$
versus $\beta$ in Example~6.2(a), (b), and (c), respectively.}

In our final example, we elaborate on the nonuniqueness
in the special case, case (iii) of Section~4, and present two distinct
sets $\{ V_1,\theta_1\}$ and $\{ V_2,\theta_2\}$
corresponding to the same scattering matrix $S.$

\noindent {\bf Example 6.3} 
As stated in Theorem~4.1(a), we note that
$\{ V_1,\theta_1\}$ and $\{ V_2,\theta_2\}$ with
$V_1(x)\equiv 0,$ $\theta_1=\pi,$ $V_1(x)\equiv 0,$ $\theta_2=\pi/2$ yield
the same scattering matrix $S_\theta(k)\equiv 1,$ as seen from
(2.8) and (2.16), illustrating the
double nonuniqueness indicated in Section~4.
We now present a less trivial example of nonuniqueness by using
the potential
$$V_1(x)=\cases 1,\qquad & 0<x<1,\\
\stretch
-a,\qquad \qquad & \ds\frac{1}{2} <x<1,\\
\stretch
0,\qquad & x>1,\endcases\tag 6.5$$
where
$a$ is a positive parameter. We can evaluate the Jost solution $f_1(k,x)$ explicitly
by using (6.5) in (2.1) and the asymptotic condition
given in (2.4) and by satisfying the continuity of
$f_1(k,x)$ and $f_1'(k,x)$ at $x=1$ and at $x=1/2.$
We then evaluate $f_1(k,0)$ and $f_1'(k,0)$ explicitly as a function of
$k$ in the presence of the parameter $a.$ Then, from (2.30) we obtain
the corresponding scattering coefficients $T_1,$ $L_1,$ and $R_1$ explicitly
via
$$T_1(k)=\ds\frac{2ik}{ik\,f_1(k,0)+f_1'(k,0)},
\quad L_1(k)=\ds\frac{ik\,f_1(k,0)-f_1'(k,0)}{ik\,f_1(k,0)+f_1'(k,0)},$$
$$R_1(k)=-\ds\frac{ik\,f_1(-k,0)+f_1'(-k,0)}{ik\,f_1(k,0)+f_1'(k,0)}.$$
We then choose the value of $a$ so that $T_1(k)$ has no poles on the positive imaginary axis
and that $T_1(0)\ne 0.$ From the small-$k$ limits of $T_1(k),$ we find that those two conditions
are satisfied provided $a$ is obtained by solving near $a=1$ the equation
$$\sqrt{a}\,\tan\left(\ds\frac{\sqrt{a}}{2}\right)=\tanh\left(\ds\frac{1}{2}\right),$$
which yields $a=0.85724\overline{7}.$ With this choice of $a,$ we get
$T_1(0)=0.97382\overline{7},$ $L_1(0)=-0.227\overline{3},$ and $R_1(0)=0.227\overline{3}.$
Note that with $a=0.85724\overline{7}$ in (6.5),
the half-line scattering matrix $S_1(k)$ corresponding to the Dirichlet
boundary condition $\theta_1=\pi$ is obtained by using the second line of (2.16).
%, and it is
%also given by the left-hand side of (4.8).
With the same specific $a$-value,
we then evaluate
the potential $V_2(x)$ corresponding to the scattering coefficients
$T_2,$ $L_2,$ $R_2,$ where
$$T_2(k)\equiv T_1(k),\quad L_2(k)\equiv -L_1(k), \quad R_2(k)\equiv -R_1(k).$$
% From (4.8) we know that the half-line scattering matrix $S_2(k)$ corresponding to the Neumann
%boundary condition $\theta_2=\pi/2$ is given by the right-hand side of (4.8). Thus,
%$S_2(k)\equiv S_1(k).$
%In fact, s
Since $T_1(k)$ has no poles on the positive imaginary axis,
one can uniquely reconstruct $V_2(x)$ from $R_2(k),$ or equivalently
 from $-R_1(k),$ with the help of (2.45), (2.44),
 and the first equation in (2.20).
Note that $V_1$ and $V_2$ can also uniquely be reconstructed from
 $L_1$ and $-L_1,$ respectively.
In fact, the corresponding
numerical approximations of $V_1$ and $V_2$ have been computed in MATLAB via the method of [18], using $L_1(k)$ and $-L_1(k)$ in the interval $k\in [0,100]$
with a discretization length of $\Delta k=0.01.$ The resulting potentials
are shown in Figure~6.2.

\vskip 10 pt

\centerline{\hbox{\psfig{figure=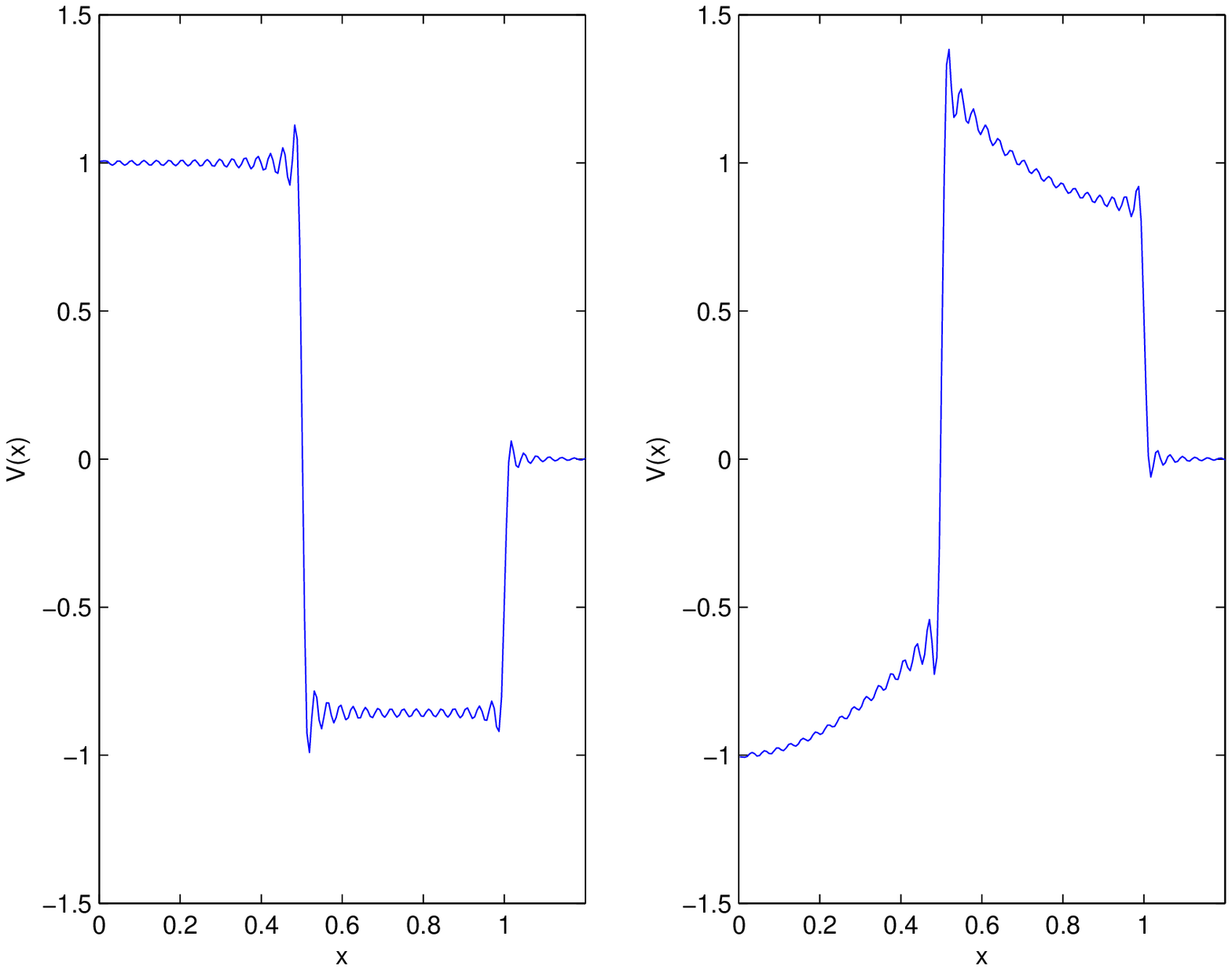,width=6.5 truein,height=3 truein}}}

{\bf Figure~6.2} The numerically reconstructed
potentials $V_1$ and $V_2$ in Example~6.3 corresponding to
$L_1$ and $-L_1,$ respectively.

\vskip 10 pt

\noindent {\bf{References}}

\vskip 3 pt

\item{[1]} Z. S. Agranovich and V. A. Marchenko,
{\it The inverse problem of scattering theory,} Gordon and Breach, New York,
1963.

\item{[2]} T. Aktosun, {\it Bound states and inverse scattering
    for
    the Schr\"odinger equation in one dimension,} J. Math. Phys.
    {\bf 35}, 6231--6236 (1994).

\item{[3]} T. Aktosun. {\it Inverse Schr\"odinger scattering on
    the line with partial knowledge of the potential,} SIAM J.
    Appl. Math. {\bf 56}, 219--231 (1996).

\item{[4]} T. Aktosun, {\it
Inverse scattering for vowel articulation with frequency-domain data,}
Inverse Problems {\bf 21}, 899--914 (2005).

\item{[5]} T. Aktosun and M. Klaus, {\it Small-energy asymptotics for the Schr\"odinger equation on the line,} Inverse
    Problems {\bf 17}, 619--632 (2001).

\item{[6]} T. Aktosun and V. G. Papanicolaou, {\it Transmission eigenvalues for the self-adjoint Schr\"odinger operator on the half
    line,} Inverse
    Problems {\bf 30}, 175001 (2014).

\item{[7]} T. Aktosun and R. Weder, {\it Inverse
    spectral-scattering problem with two sets of discrete
    spectra for the radial Schr\"odinger equation,} Inverse
    Problems {\bf 22}, 89--114 (2006).

\item{[8]} K. Chadan and P. C. Sabatier, {\it Inverse problems
    in quantum scattering theory,} 2nd ed., Springer, New York,
    1989.

\item{[9]} P. Deift and E. Trubowitz, {\it Inverse scattering
    on the line,} Commun. Pure Appl. Math. {\bf 32}, 121--251
    (1979).

\item{[10]} L. D. Faddeev, {\it Properties of the $S$-matrix of
    the one-dimensional Schr\"odinger equation,} Amer. Math.
    Soc. Transl. {\bf 65} (ser. 2), 139--166 (1967).

\item{[11]} I. M. Gel'fand and B. M. Levitan, {\it On the
    determination of a differential equation from its spectral
    function,} Amer. Math. Soc. Transl. {\bf 1} (ser. 2),
    253--304 (1955).

\item{[12]} G. M. L. Gladwell, {\it Inverse problems in vibration,} 2nd ed.,
Kluwer, Dordrecht, 2004.

\item{[13]} B. Gr\'ebert and R. Weder,
{\it Reconstruction of a potential on the line
that is a priori known on the half line,}
SIAM J. Appl. Math. {\bf 55}, 242--254 (1995).

\item{[14]} B. M. Levitan, {\it Inverse Sturm-Liouville
    problems,} VNU Science Press, Utrecht, 1987.

\item{[15]} V. A. Marchenko, {\it Sturm-Liouville operators and
    applications,} Birkh\"auser, Basel, 1986.

\item{[16]} R. G. Newton, {\it Scattering theory of waves and
    particles,}
    2nd ed., Springer, New York, 1982.

    \item{[17]} N. N. Novikova and V. M. Markushevich, {\it
    Uniqueness of the solution of the
    one-dimensional problem of scattering for potentials
    located on the positive semiaxis,} Comput. Seismol.
{\bf 18}, 164--172 (1987).

\item{[18]} P. E. Sacks, {\it Reconstruction of steplike potentials,}
Wave Motion {\bf 18}, 21--30 (1993).

\item{[19]} M. Zworski,
{\it Distribution of poles for scattering on the real line,}
J. Func. Anal. {\bf 73}, 277--296 (1987).

\end